\documentclass[a4paper,10pt]{article}
\usepackage{amsmath, amssymb}

\usepackage[francais,english]{babel}

\usepackage{amsmath, amssymb}
\usepackage{geometry,amsmath,amssymb,enumerate,amsthm}

\usepackage{amsfonts}
\usepackage{amsmath,amsthm,indentfirst}
\usepackage{amssymb}
\usepackage{amsthm}
\usepackage{amssymb}
\usepackage{amsmath}
\usepackage{fancyhdr}

\usepackage[pdftex,bookmarks,colorlinks,breaklinks]{hyperref}

\newcommand{\w}{w_{\beta}}

\newcommand{\g}{\gamma_{\beta}}
\newcommand{\dw}{{w\hspace{0,0098889cm}dx}}

\newtheorem{theorem}{Theorem}[section]

\newtheorem{lemma}[theorem]{Lemma}

\theoremstyle{definition}

\newtheorem{remark}[theorem]{Remark}

\numberwithin{equation}{section}

\begin{document}

\date{}
\author{Omar Lazar and Pierre-Gilles Lemari\'e-Rieusset  
 }
\title{Infinite energy solutions for a 1D transport equation with nonlocal velocity } 
\maketitle
\bibliographystyle{plain}

\begin{abstract}
We study a one dimensional dissipative transport equation with  nonlocal velocity and critical dissipation. We consider the Cauchy problem for initial values with infinite energy. The control we shall use involves 
  some weighted Lebesgue  or Sobolev spaces.  More precisely,  we consider the family of weights   given by $w_{\beta}(x)=(1+\vert x \vert^{2})^{-\beta/2}$ where $\beta$ is a real parameter in $(0,1)$ and we treat the Cauchy problem for the cases    $\theta_{0} \in H^{1/2} (w_{\beta})$ and $\theta_{0} \in H^{1} (w_{\beta})$ for which we prove global existence results (under smallness assumptions on the $L^\infty$ norm of $\theta_0$). The key step in the proof  of our theorems is based on the use of  two new  commutator estimates involving fractional differential operators and the family of  Muckenhoupt weights. 
\end{abstract}

\maketitle
\tableofcontents

\section*{Introduction}
In this article, we are interested in the following 1D transport equation with nonlocal velocity which has  been introduced by C\'ordoba, C\'ordoba and Fontelos in \cite{CCF} :
\begin{equation} 
\ (\mathcal{T}_{\alpha}) \ : \\\left\{
\aligned
&\partial_{t}\theta+\theta_x \mathcal{H}\theta+ \nu \Lambda^{\alpha}\theta = 0 \hspace{2cm} 
\\ \nonumber
& \theta(0,x)=\theta_{0}(x).
\endaligned
\right.
\end{equation}

\noindent Here, $\mathcal{H}$ denotes the Hilbert transform,  defined by 
\begin{equation} \label{hil}
\mathcal{H} \theta \equiv \frac{1}{\pi} PV \int \frac{\theta(y)}{{y-x} } \ dy,
\end{equation}
and the operator $\Lambda^{\alpha}$ is defined (in 1D) as follows
\[\Lambda^{\alpha}\theta \equiv (-\Delta)^{\alpha/2} \theta=C_{\alpha}P.V.\int_{\mathbb{R}}{\frac{\theta(x)-\theta(x-y)}{|y|^{1+\alpha}}dy} \]
where $C_{\alpha}>0$ is a positive constant  which depends only on $\alpha$ and $0<\alpha< 2$ is a real parameter. Note that with this convention in \ref{hil}, we have $\partial_{x} \mathcal{H}=-\Lambda$ \\

This equation can be viewed as a toy model for several equations coming from  problems in fluid dynamics, in particular it models the 3D Euler equation written in vorticity form (see e.g. \cite{CLM}, \cite{Bak}, \cite{De1}, \cite{OSW} where other 1D models for 3D Euler equation are studied).

 One can observe that this equation is a one dimensional model for the 2D dissipative Surface-Quasi-Gesotrophic $({SQG})_{\alpha}$ equation (see \cite{CMT}) written in a non-divergence form (see also \cite{CaC}, \cite{CaC2}, \cite{ChCCF}, \cite{RL2}  where the divergence form equation is studied).  The 2D dissipative $SQG$ equation reads as follows
\begin{equation}    
({SQG})_{\alpha} : \ \\\left\{
\aligned
&\partial_{t}\theta(x,t)+ u( \theta) .\nabla\theta+ \nu\Lambda^{\alpha}\theta = 0
\\ \nonumber
& \theta(0,x)=\theta_{0}(x),
\endaligned
\right.
\end{equation}  
where the velocity $u(\theta)= \mathcal{R}^{\perp}\theta$ is given by the Riesz transforms $\mathcal{R}_1\theta$ and $\mathcal{R}_2\theta$ of $\theta$  as  $$u(\theta)=(-\mathcal{R}_2\theta,\mathcal{R}_1\theta)= (-\partial_{x_2}
\Lambda^{-1}\theta,
\partial_{x_1}\Lambda^{-1}\theta).$$ Obviously the velocity $u(\theta)$  is divergence free. In 1D, we lose this divergence free condition, while the analogue of the Riesz transforms is the Hilbert transform;    one gets the equation $(\mathcal{T}_{\alpha})$. \\
   
 One can also see this equation as an analogue of the fractional Burgers equation with the nonlocal velocity $u(\theta)=\mathcal{H} \theta$ instead of $u(\theta)= \theta$.  However, the nonlocal character of the velocity makes the $(\mathcal{T}_{\alpha})$ equation  more complicated to deal with comparing to the fractional Burgers equations which is now quite well understood (see \cite{KNS},  \cite{CCS}, \cite{Ki}). Finally, let us mention that this equation also shares some similarities with the Birkhoff-Rott equation which modelises the evolution of a vortex patch, we refer to \cite{CCF}, \cite{Bak} for more details regarding this analogy. \\

It is easy to guess that this 
  kind of fractional transport equation admits  an $L^{\infty}$ maximum principle (due to the diffusive character of $-\Lambda^\alpha$ and the presence of the derivative $\theta_x$ in the advection term).  For $\theta \in  L^\infty$, one thus may view $\theta_x \mathcal{H} \theta$ as a term of order $1$, while $\Lambda^\alpha$ is of order $\alpha$; thus,  one has to consider 3 cases depending on the value of $\alpha$,  namely  $\alpha \in (0,1)$, $\alpha=1$ and $\alpha \in (1,2)$. They are respectively called supercritical, critical and subcritical cases. \\
  
   The inviscid case  (i.e.  $\nu=0$) was first studied by C\'ordoba, C\'ordoba and Fontelos in \cite{CCF} where the authors proved that blow-up of regular solutions may occur. They  proved that there exists a family of smooth, compactly supported, even and positive initial data for which the associated solution blows up in finite time.  By adapting the method used in \cite{CCF} along with the use of new nonlocal inequalities obtained in \cite{CCF2}, Li and Rodrigo \cite{RL}  proved that blow-up of smooth solutions also holds in the viscous case, in the range $\alpha \in (0,1/2)$. Using a different method,  Kiselev \cite{Ki} was able to prove that singularities may appear in the case $\alpha \in [0,1/2)$ (where the case $\alpha=0$ conventionnally designs the inviscid case $\nu=0$).  In this latter range, that is $\alpha \in [0,1/2)$,  Silvestre and Vicol \cite{ViS} gave four differents proofs of the same results as \cite{CCF}, \cite{RL}, \cite{ViS},  namely they proved the existence of singularities for classical ($\mathcal{C}^{1}$)   solutions starting from a well chosen class of initial data. In \cite{Do}, T. Do showed eventual regularization in the supercritical case and global regularity for the slightly supercritical version of equation $\mathcal{T}_{\alpha}$,   in the spirit of what   was done for the SQG equation in \cite{S}, \cite{Ki}. One can also see the articles \cite{Dong} and \cite{BG} where local existence results are obtained in this regime. In the range $\alpha \in [1/2,1)$, the question about blow-up or global existence of regular solutions remains open. 
   
    The critical and the sub-critical cases are well understood. Indeed, by adapting methods introduced in \cite{KNV}, \cite{CV}, \cite{CVi},  one recovers all the results known for the critical SQG equation, under an extra  positiveness assumption on the initial data (see \cite{K}). The first global existence results are those of C\'ordoba, C\'ordoba and Fontelos \cite{CCF}. They obtained global existence results for non-negative  data in $H^{1}$ and $H^{1/2}$  in the subcritical case and also in the critical case under a smallness assumption of the $L^{\infty}$ norm of the initial data. In \cite{Dong}, Dong treated the critical case and obtained the global well-posedness  for data in  $H^{s}$  where $s >3/2-\alpha$ and without sign conditions on the  initial data. In the critical case, Kiselev proved in \cite{Ki} that  there exists a unique global smooth solution for all $\theta_{0} \in H^{1/2}$.   \\
    
  In this article, we will  focus on the critical case ($\alpha=1$) and without loss of generality we shall fix $\nu=1$. Futhermore, in contrast with \cite{CCF} and \cite{Dong}, we shall not assume that $\theta$ decays at infinity fast enough to ensure that $\|\theta\|_2<+\infty$. It is worth pointing out that, our solutions being of infinite energy, one cannot directly use methods coming from $L^{\infty}$-critical case used for instance in  \cite{CV}. However, in the case of an infinite-energy data, one can still use energy estimates (in the spirit of \cite{CCF})  to prove global existence results provided that $\theta$ decreases only at a slow rate, namely
$$
\int \vert \theta(x,t) \vert^{2} \frac{dx}{(1+\vert x \vert^{2})^{\beta/2}}  <+\infty
 $$
 The weight we consider is therefore given by $w_{\beta}(x)=(1+\vert x \vert^{2})^{-\beta/2}$.  Motivated by the work done in \cite{CCF}, we  will   study the cases of   small data in  $L^\infty$ which belong moreover to    $H^{1/2}(w_{\beta})$ or  $H^{1}(w_{\beta})$, although one can generalize to a higher regularity class of initial data (we think that it should be even easier to treat).  When the initial data lies in $H^{1/2}(w_{\beta})$ or  $H^{1}(w_{\beta})$ we prove global existence  of weighted Leray-Hopf type  solutions but we require the $L^{\infty}$ norm of the initial data to be small enough.  As one may expect, in the subcritical case one can prove the existence of global solutions without smallness assumption. This is done by the first author in \cite{Laz} using Littlewood-Paley theory along with a suitable commutators estimates.  He also treated the supercritical case where he obtained local existence results for arbitrary big initial data \cite{Laz}. \\
   
 The construction of a solution is  based on an energy method and amounts to control some nontrivial commutators involving the weight $w_{\beta}$ along with some classical harmonic analysis tools such as the use of the Hardy-Littlewood maximal function and Hedberg's inequality for instance (see \cite{Hed}, \cite{Stein});  such tools are motivated by the fact $w_{\beta}$ is a Muckenhoupt weight. The new commutator estimates can be used to prove existence of infinite energy solutions for other nonlinear transport equations with fractional diffusion such as the 2D dissipative quasi-geostrophic equation as well as  the fractional porous media equation for instance. \\
 
  The rest of the paper is organized into five sections.  In the first section, we  state our main theorems. In the second section we recall some results concerning the Muckenhoupt weights. In the third and fourth section, we respectively establish {\it{a priori}} estimates and prove our main results. In the last section we revisit the construction of regular enough solutions.
    
 \section{Main theorems}
 
 In the case of a weighted  $H^{1/2}$ data we have the following theorem,
 \begin{theorem} \label{th22}
Let $0<\beta<1$ and $\w(x)=(1+x^2)^{-\beta/2}$. There exists a constant $C_{\beta}>0$ such that, whenever $\theta_{0}$ satisfies the conditions 
\begin{itemize}
\item $\theta_{0}$ is bounded and small enough : $\vert \theta_0\vert  \leq C_{\beta}$ 
\item $\displaystyle\int \vert\theta_{0}\vert^{2} \w(x) \  dx < \infty$  and $\displaystyle\int \vert \Lambda^{1/2} \theta_{0}\vert^{2} w_\beta(x) \  dx < \infty$, \\
\end{itemize} 
 there exists a solution $\theta$ to equation $\mathcal{T}_1$ such that, for every $T>0$, we have 

\begin{itemize}
\item $\displaystyle\sup_{0<t<T} \int \vert\theta(t,x)\vert^{2} \w(x) \  dx < \infty$ 
\item $\displaystyle\sup_{0<t<T} \int \vert\Lambda^{1/2} \theta(t,x)\vert^{2} \w(x) \  dx < \infty$ 
\item $\displaystyle\int_{0}^{T} \int \vert \Lambda \theta(t,x) \vert^{2} \w(x) \ dx \ dt < \infty$ \\
\end{itemize} 

\end{theorem}

 A similar result holds for higher regularity (weighted $H^1$ data).
 
  \begin{theorem} \label{th23}
Let $0<\beta<1$ and $\w(x)=(1+x^2)^{-\beta/2}$. There exists $C_{\beta}>0$ such that, whenever $\theta_{0}$ satisfies the conditions 
\begin{itemize}
\item $\theta_{0}$ is bounded and small enough : $\vert \theta_0\vert \leq C_{\beta}$ 
\item $\displaystyle\int \vert\theta_{0}\vert^{2} \w(x) \  dx < \infty$  and $\displaystyle\int \vert \Lambda \theta_{0}\vert^{2} \w(x) \  dx < \infty$, \\
\end{itemize} 
  there exists a solution $\theta$ to equation $\mathcal{T}_1$ such that, for every $T>0$, we have 

\begin{itemize}
\item $\displaystyle \sup_{0<t<T} \int \vert\theta(t,x)\vert^{2} \w(x) \  dx < \infty$ 
\item $\displaystyle \sup_{0<t<T} \int \vert\Lambda \theta(t,x)\vert^{2} \w(x) \  dx < \infty$ 
\item $\displaystyle \int_{0}^{T} \int \vert \Lambda^{3/2} \theta(t,x) \vert^{2} \w(x) \ dx \ dt < \infty$ \\
\end{itemize} 

\end{theorem}

\section{Preliminaries on the Muckenhoupt weights.}

In this section, we briefly recall the tools and the notations we shall use throughout the article. We first recall some basic facts and notations on weighted Lebesgue or Sobolev spaces. A weight $w$ is a positive and locally integrable function. A measurable function $\theta$ on $\mathbb R$ belongs to the weighted Lebesgue spaces $L^p(w dx)$ with $1\leq p < \infty$ if and only if
$$
\Vert \theta \Vert_{L^{p}(\dw)} = \left( \int \vert \theta (x) \vert^{p} \ w(x) \ dx \right)^{1/p} < \infty.
$$
An important class of weights is the so-called Muckenhoupt class $\mathcal{A}_{p}$ for $1<p<\infty$. A weight is said to be in the $\mathcal{A}_{p}$ class of Muckenhoupt (with $p \in (1,\infty)$) if and only if there exists a constant $C(w,p)$ such that we have the reverse H\"older inequality
$$
\sup_{r>0, x_{0} \in \mathbb R} \left( \frac{1}{2r} \int_{[x_{0}-r, x_{0}+r]} w(x) \ dx \right) \left( \frac{1}{2r} \int_{[x_{0}-r, x_{0}+r]} w(x)^{-\frac{1}{p-1}} \ dx \right)^{p-1} \leq C(w,p).
$$In particular, if $0<\beta<1$, then the weight $w_{\beta}(x)=(1+\vert x \vert^{2})^{-\beta/2}$ belongs to the $\mathcal{A}_{p}$ class for all $1<p<\infty$.

  Let us recall that the Hardy-Littlewood maximal function of a locally integrable function $f$ on $\mathbb R$ is defined by
$$
\mathcal{M} f (x)= \sup_{r>0} \frac{1}{2r} \int_{[x-r,x+r]} \vert f(y) \vert \ dy.
$$
We have the following other characterization of the $\mathcal{A}_{p}$ class \cite{CM2}, \cite{Muc} : a  weight $w$ belongs to $  \mathcal{A}_{p}$ if and only if there exists a constant $C_{p,w}$ such that for every $f \in L^{p}(w \ dx)$, we have
$$
\Vert \mathcal{M} f (x) \Vert_{L^{p}(w \ dx)} \leq C_{p,w} \Vert f \Vert_{L^{p}(\dw)}.
$$

\noindent Another important property of Muckenhoupt weights is that Calder\'on-Zygmund type operators are bounded on $L^{p}(w \ dx)$ when $w \in \mathcal{A}_{p}$ and $1<p<\infty$. We shall use this property  in the case of the Hilbert transform $\mathcal{H}$  and in the case of the truncated Hilbert transform, defined by
\begin{equation} \label{bh}
\mathcal{H}_{\#} f(x) = \frac{1}{\pi} P. V. \int \frac{\alpha(x-y)}{x-y} f(y) \ dy
\end{equation}
where $\alpha$ is an  even, smooth and compactly supported function such that $\alpha(x)=1$ if $\vert x \vert<1$ and $\alpha(x)=0$ if $\vert x \vert>2$. We refer for instance to \cite{Stein} or \cite{HMW} for more details.\\

We now recall the definition of the weighted Sobolev spaces $H^{1}( \dw)$ and $H^{1/2}(\dw )$. The space $H^{1}(w \ dx)$ is defined by
$$
f \in H^{1}( \dw) \Leftrightarrow f \in L^{2}(\dw) \ \  {\rm{and}} \ \ \partial_{x} f \in L^{2}(\dw).
$$
Note that, due to \ref{hil},  we have
$$
\mathcal{H} \partial_x = \Lambda \ \  {\rm{and}} \ \  \mathcal{H}\Lambda = \partial_x,
$$
we see that, when $w \in \mathcal{A}_{2}$, the semi-norm $\Vert \partial_{x} f \Vert_{L^{2}(w \ dx)}$ is equivalent to the semi-norm  $\Vert  \Lambda  f \Vert_{L^{2}( \dw)}$. Therefore, when $w \in \mathcal{A}_{2}$, we have the following equivalence
\begin{equation*} \label{h1}
f \in H^{1}( \dw) \Leftrightarrow (1-\partial^{2}_{x})^{1/2} f \in L^{2}(  \dw) \Leftrightarrow f \in L^{2}(\dw) \ \  {\rm{and}} \ \ \Lambda f \in L^{2}(\dw).
\end{equation*}
Analogously, we define the spaces $H^{1/2}(\dw )$ as
\begin{equation*}
f \in H^{1/2}(w \ dx) \Leftrightarrow (1-\partial^{2}_{x})^{1/4} f \in L^{2}(\dw) \Leftrightarrow f \in L^{2}(\dw) \ \  {\rm{and}} \ \ \Lambda^{1/2} f \in L^{2}(\dw).
\end{equation*}

The following useful property will be used several times (see \cite{Stein}, p.57). Fix an integrable nonnegative and radially decreasing function $\phi$ such that its integral over $\mathbb R$ is equal to 1. We set,  $\phi_{k}(x)=k^{-1} \phi(xk^{-1})$ for all $k>0$, then 
\begin{equation} \label{St}
\sup_{k>0} \vert f * \phi_{k}(x) \vert \leq \mathcal{M}f(x).
\end{equation}
In the sequel, we shall use Gagliardo-Nirenberg type inequalities in the weighted setting. Let us first note that, provided $f$ vanishes at infinity (in the sense that $\displaystyle\lim_{t\rightarrow +\infty} e^{t\Delta} f=0$ in $\mathcal{S}'$), one may write
\begin{equation*}
-f = \int_{0}^{\infty} e^{t\Delta} \Delta f \ dt.
\end{equation*}
where $f \mapsto e^{t\Delta}f$ is the heat kernel operator defined by $e^{t\Delta}f=G(x,t)*f$ where $*$ is the convolution with respect to the $x$ variable and $G(x,t)=(4\pi t)^{-1/2} e^{-\frac{x^2}{4t}}$ which verifies the heat equation $\partial_{t} G(x,t)=\Delta G(x,t)$.  \\

\noindent Then, for all $N \in \mathbb{N}^*$ by writing $1=\partial^{N-1}_{t} (\frac{t^{N-1}}{(N-1)!})$ and  integrating by parts $(N-1)$ times, one obtain the following equality
$$
-f= \frac{1}{(N-1)!} \int_{0}^{\infty} (-t \Delta)^{N} e^{t \Delta} f \frac{dt}{t}.
$$ 
Then,  for $0<\gamma < \delta < 2N$,  using the fact that the operator $\Lambda^{2N-\delta+\gamma}$ is a convolution operator with an integrable kernel which is dominated by an integrable radially decreasing function,  along with the inequality
$$
\sup_{t} \vert \Lambda^{\gamma} e^{t\Delta} f \vert \leq c t^{-\gamma/2} \mathcal{M} f(x)
$$
allow us to get
$$
\vert \Lambda^{\gamma} f(x) \vert \leq C \int_{0}^{\infty} \min(t^{-\gamma/2} \Vert f \Vert_{\infty}, t^{\frac{\delta-\gamma}{2}} \mathcal{M}(\Lambda^{\delta}f)(x)) \frac{dt}{t}
$$
Then, we recover Hedberg's inequality (see Hedberg \cite{Hed})
\begin{equation} \label{Hed}
\vert \Lambda^{\gamma} f(x) \vert \leq C_{\gamma, \delta} ( \mathcal{M}(\Lambda^{\delta}f)(x)))^{\frac{\gamma}{\delta}} \Vert f \Vert^{1-\frac{\gamma}{\delta}}_{\infty}
\end{equation}
Note that, if $\gamma \in \mathbb{N}^{*}$, one may replace $\Lambda^{\gamma} f(x)$ with $\partial^{\gamma}_{x} f(x)$. Using  (\ref{Hed}), one easily deduce  the following Gagliardo-Nirenberg type inequalities provided that the weight $w \in \mathcal{A}_3$ (actually $w \in \mathcal{A}_2$ suffices for \ref{b1} and \ref{b})  \
\begin{eqnarray} \label{GN}
&& \Vert \Lambda^{1/2} f \Vert_{L^{4}(w  dx)} \leq C \Vert f \Vert^{1/2}_{\infty} \Vert \Lambda f \Vert^{1/2}_{L^{2}(w dx)} \label{b1} \\
&& \nonumber \\ 
&& \Vert \Lambda f \Vert_{L^{3}(w  dx)} \leq C \Vert f \Vert^{1/3}_{\infty} \Vert \Lambda^{3/2} f \Vert^{2/3}_{L^{2}(w dx)} \label{b} \\
&&\nonumber \hspace{-4cm} {\rm{and}} \\
&& \Vert \partial_x f \Vert_{L^{3}(w  dx)} \leq C \Vert f \Vert^{1/3}_{\infty} \Vert \Lambda^{3/2} f \Vert^{2/3}_{L^{2}(w dx)}\label{b2}
\end{eqnarray}\\

\noindent For instance, to prove \ref{b}, it suffices to set $\gamma=1$ and $\delta=3/2$ in  \ref{Hed}  and to raise to the power 3 in both side, one obtain 
$$
\vert \Lambda f\vert^{3} \leq ( \mathcal{M}(\Lambda^{3/2}f)(x))^2 \Vert f \Vert_{\infty}.
$$
Mutiplying this latter inequality by $w$ and integrating with respect to $x$ give
$$
\Vert \Lambda f \Vert^{3}_{L^{3}(w dx)} \leq \Vert f \Vert_{\infty} \Vert \mathcal{M}(\Lambda^{3/2}f)(x)\Vert^{2}_{L^{2}(wdx)} \leq  \Vert f \Vert_{\infty} \Vert \Lambda^{3/2}f\Vert^{2}_{L^{2}({wdx})},
$$
where, in the last inequality, we used  the continuity of the maximal function $\mathcal{M}$ on $L^{2}(wdx)$ because $w\in \mathcal{A}_2$. Therefore, inequality \ref{b} follows by taking the power $1/3$ in both sides. Then,  observe that \ref{b2} is a direct consequence of \ref{b}. Indeed, we have   $\mathcal{H}\Lambda  f=-\partial_{x}f$ and due to  $w\in \mathcal{A}_{3}$ one can use the continuity on $L^{3}({wdx})$ of $\mathcal{H}$ to obtain the inequality
$$
\Vert \partial_{x}f\Vert_{L^{3}(wdx)}=\Vert \mathcal{H}\Lambda  f \Vert_{L^{3}(wdx)}\leq\Vert \Lambda  f \Vert_{L^{3}(wdx)},
$$
therefore we recover \ref{b2}. \\

The space of positive smooth functions compactly supported in an open set $\Omega$ will be denoted by $\mathcal{D}(\Omega)$. We shall use the notation $A \lesssim B$  if there exists  constant $C>0$ depending only on controlled quantities such that $A \leq C B$. We shall often use the same notation to design a controlled constant although it is not the same from a line to another. Note that we shall write indifferently $\partial_x\theta$ or $\theta_x$ for the derivative  as well as $\Vert . \Vert_{p}$ or $\Vert . \Vert_{L^{p}}$ for the classical Lebesgue spaces.  \\

\section{Useful lemmas}

In our future estimations, we will need to control the $L^{p}$ norm of some nontrivial commutators involving our weight $w_{\beta}$ and the nonlocal operators $\Lambda$ and $\Lambda^{1/2}$. Also, a control of $\Lambda w$ by $c w$ will be needed. The aim of this section is to establish all those nonlocal estimates involving $w$ and $\Lambda$. Before starting the proofs of those commutator estimates, we shall give some important remarks that will be helpful  to estimate singular integrals involving the weight $w$. When estimating commutators involving the weight $\w(x)=(1+x^2)^{-\beta/2}$, we are lead to estimate quantities such that $\w(x)-\w(y)$.  In order to estimate  $\w(x)-\w(y)$, we shall distinguish three areas that we will call $\Delta_{1}(x), \Delta_{2}(x)$ and $\Delta_{3}(x)$. Those areas are defined as follows,
\begin{eqnarray*}
 && \Delta_{1}(x) = \{ y \  /  \ \vert x - y \vert < 2 \} \\
&& \Delta_{2}(x) = \{ y \  /  \ \vert x - y \vert \geq 2 \}  \cap \{  y \  /   \vert x - y \vert \leq \frac{1}{2} \max(\vert x \vert, \vert y \vert) \} \\  
&& \Delta_{3}(x) =  \{ y \  /  \ \vert x - y \vert \geq 2 \}  \cap \{  y \  /   \vert x - y \vert > \frac{1}{2} \max(\vert x \vert, \vert y \vert) \}.
 \end{eqnarray*}
 Note that we have $\mathbb R =  \Delta_{1}(x) \cup \Delta_{2}(x)  \cup \Delta_{3}(x) $. In the sequel, we shall also use the notation $\w(x)\approx \w(y)$ if there exists two positive  constants $c$ and $C$ such that $c \leq \frac{w(x)}{w(y)} \leq C$. In those different areas, we will need to use the following estimates :
 \begin{itemize}
 \item A straightforward computation gives that   $\vert \partial_{x} w_{\beta}(x) \vert +\vert\partial^2_x w_\beta(x)\vert \leq C w_{\beta}(x)$
 \item On $\Delta_{1}(x)$, we have that $\w(x)\approx \w(y)$ and moreover
  $$\vert \w(x)-\w(y) \vert \leq \vert x-y\vert\sup_{z\in[x,y]}\vert\partial_xw_\beta(z)\vert\leq  C \vert x -y \vert w_\beta(x)$$
  On the other hand, if $\alpha$ is an even, smooth and compactly supported function such that $\alpha(x)=1$ if $\vert x \vert <1$ and $\alpha(x)=0$ if $\vert x \vert >2$, then
\begin{equation} \label{trunc}
 \vert \w(y)-\w(x) + \alpha(x-y)(x-y) \partial_{x} \w(x) \vert \leq C \vert x-y \vert^{2} \w(x)
\end{equation} 
 \item On $\Delta_{2}(x)$, we shall only use that  $\w(x)\approx \w(y)$ 
 \item On $\Delta_{3}(x)$, we have $1 \leq \w(x)^{-1} \leq C \vert x-y \vert^{\beta}$ and $1 \leq \w(y)^{-1} \leq C \vert x-y \vert^{\beta}$
  \end{itemize}
  
  \begin{remark}
  Obviously, similar estimates hold for $\g(x)=\w(x)^{1/2}$. Indeed, it suffices to replace $\w$ with $\g$ and $\beta$ with $\beta/2$. 
  \end{remark}
  
 The purpose of the following subsections is to prove that we can indeed control those commutators and that we have a nice bound for $\Lambda w_{\beta}$.

\subsection{Two commutator estimates involving the weight $w_{\beta}$}

In the next lemma, we obtain two  commutator estimates that are crucial in the proof of the energy inequality. 

\begin{lemma} \label{comm} Let $w_{\beta}(x)=(1+x^2)^{-\beta/2}$, $0<\beta<1$, then we have the two following estimates
\begin{itemize}
\item Let $p\geq2$ be such that $\frac{3}{2}-\beta(1-\frac{1}{p}) >1$, then the commutator $\frac{1}{w_{\beta}} [\Lambda^{1/2},w_{\beta}]$ is bounded from $L^{p}(\w dx)$ to $L^{p}(\w dx)$.
 \item Let $2\leq p < \infty$,  then the commutator $\frac{1}{\sqrt{\w}} [\Lambda,\sqrt{\w}]$ is bounded from $L^{p}(\w dx)$ to $L^{p}(\w dx)$.
 \end{itemize}
 \end{lemma}

  \noindent ${\bf{Proof \ of \ lemma \ \ref{comm}}.}$ \\

\noindent Let us prove the first commutator estimate. We first write
$$
\Lambda^{1/2} f(x)=c_0 \int \frac{f(x)-f(y)}{\vert x-y \vert^{3/2}} \ dy
$$
 so that
$$
\frac{1}{\w(x)} [\Lambda^{1/2}, \w] f(x) = c_0 \frac{1}{\w(x)^{\frac{1}{p}}} \int \frac{\w(x)-\w(y)}{\w(x)^{1-\frac{1}{p}} \w(y)^{\frac{1}{p}} \vert x-y \vert^{3/2} } \w(y)^{\frac{1}{p}} f(y) \ dy
$$
Let us set
$$
K(x,y) \equiv  \frac{\w(x)-\w(y)}{\w(x)^{1-\frac{1}{p}} \w(y)^{\frac{1}{p}} \vert x-y \vert^{3/2} }
$$
On $\Delta_{1}(x)$ we have 
$$ 
\vert K(x,y) \vert \leq C \frac{1}{\vert x-y \vert^{1/2}}
$$
On $\Delta_{2}(x)$, since $\w(x) \approx \w(y)$, we get 
$$
\vert K(x,y) \vert \leq C \frac{1}{\vert x-y \vert^{3/2}}
$$
On $\Delta_{3}(x)$, we have the following estimate
$$
\vert K(x,y) \vert \leq C \frac{\w(x)^{\frac{1}{p}-1} + \w(y)^{-\frac{1}{p}}}{\vert x-y \vert^{3/2}} \leq C' \frac{1}{\vert x-y \vert^{\frac{3}{2}-\beta(1-\frac{1}{p})}}
$$
Note that, for $0<\beta<1$ we have $\frac{3}{2}-\beta(1-\frac{1}{p})>1$ if $p\geq2$. Therefore, if we introduce the function $x \mapsto \Phi(x)$ as follows
$$
\Phi(x) \equiv \min\left(\frac{1}{\vert x \vert^{1/2}},{\frac{1}{\vert x \vert^{\frac{3}{2}-\beta(1-\frac{1}{p})}}}\right),
$$
we find that $\Phi$ belongs  to $L^{1}(\mathbb R)$ and that
$$
\left \vert \frac{1}{\w(x)} [\Lambda^{1/2}, \w] f(x) \right \vert \leq C \frac{1}{\w(x)^{\frac{1}{p}}} \int \Phi(x-y) \w(y)^{\frac{1}{p}} \vert f(y) \vert \ dy
$$
The integral appearing in the right hand side is nothing but the convolution of $x\mapsto \Phi(x) \in L^{1}(\mathbb R)$ with $x\mapsto \w(x)^{\frac{1}{p}} \vert f(x) \vert \in L^{p}(\mathbb R)$. To finish the proof, we just have to take the power $p$ in both side  then to integrate with respect to $x$ and by Young's inequality for convolution, we get
$$
\int \left \vert \frac{1}{\w(x)} [\Lambda^{1/2}, \w] f(x) \right \vert^{p} \ w  dx \leq C \int \left\vert (\Phi * \w^{1/p} f)(x) \right\vert^{p} \ dx \leq C \Vert \Phi \Vert^{p}_{L^{1}} \Vert \w^{1/p} f \Vert^{p}_{L^{p}}
$$
and therefore, 
$$
\left\Vert \frac{1}{\w(x)} [\Lambda^{1/2}, \w] f(x) \right \Vert_{L^{p}(\w dx)} \leq C \Vert f \Vert_{L^p(\w dx)}
$$ 
Let us prove the second commutator estimate. Let us denote $\g=\sqrt{\w}$, recall that 
$$
\Lambda f (x)=\frac{1}{\pi} \lim_{\epsilon \rightarrow 0} \int_{\epsilon < \vert x-y \vert <\frac{1}{\epsilon}} \frac{f(x)-f(y)}{\vert x-y \vert^{2}} \ dy
$$
Therefore,
$$
\frac{1}{\sqrt{\w(x)}} [\Lambda, \sqrt{\w(x)} ] f= \frac{1}{\pi \g(x)} \lim_{\epsilon \rightarrow 0} \int_{\epsilon < \vert x-y \vert <\frac{1}{\epsilon}} \frac{\g(y)-\g(x)}{\vert x-y \vert^{2}} \ f(y)\ dy  
$$
Then, as we did before, we split the integral into three pieces. In other to deal with the integration in $\Delta_{1}(x)$,  we need to introduce a even, smooth and compactly supported function $\alpha$ such that $\alpha(x)=1$ if $\vert x \vert <1$ and $\alpha(x)=0$ if $\vert x \vert >2$. By doing so, we get an extra term which is nothing but the truncated Hilbert transform  of $f$ (see \ref{bh}) times another controlled term.   More precisely, we  write the commutators as follows
\begin{eqnarray*}
\frac{1}{\sqrt{\w(x)}} \left[\Lambda, \sqrt{\w(x)} \right] f &=&  \frac{1}{\pi \g(x)} \int_{\Delta_{1}(x)} \frac{\g(y)-\g(x)-(y-x)\alpha(x-y) \partial_{x} \g(x)}{\vert x-y \vert^{2}} f(y) \ dy \\
 &-& \frac{\partial_{x} \g(x)}{\g(x)} \mathcal{H}_{\#} f(x) +  \frac{1}{\pi \g(x)} \int_{\Delta_{2}(x) \cup \Delta_{3}(x)} \frac{\g(y)-\g(x)}{\vert x-y \vert^{2}} \ f(y) \ dy
\end{eqnarray*}
Then, observe that on $\Delta_{1}(x)$ we have (see \ref{trunc})
$$
\frac{1}{\g(x)^{1-\frac{2}{p}} \g(y)^{\frac{2}{p}}} \frac{\vert \w(y)-\w(x) + \alpha(x-y)(x-y) \partial_{x} \w(x) \vert}{\vert x-y \vert^{2}} \leq C
$$
On $\Delta_{2}(x)$, we have
$$
\frac{1}{\g(x)^{1-\frac{2}{p}} \g(y)^{\frac{2}{p}}} \frac{\vert \g(y)-\g(x)\vert}{\vert x-y\vert^{2}} \leq C \frac{1}{\vert x-y \vert^{2}}
$$

\noindent  Here, we used the property that on $\Delta_{1}(x)$ and $\Delta_{2}(x)$ we have $\g(x) \approx \g(y)$ and therefore $\g(x)^{1-\frac{2}{p}} \g(y)^{\frac{2}{p}} \approx \g(x)$. \\

\noindent Finally, on $\Delta_{3}(x)$ we use the fact that $\g(x) \leq 1, \g(x)^{-1} \leq C \vert x -y \vert^{\beta/2}$. We also have that   $\g(y) \leq 1, \g(y)^{-1} \leq C \vert x -y \vert^{\beta/2}$, therefore
\begin{eqnarray*}
\frac{1}{\g(x)^{1-\frac{2}{p}} \g(y)^{\frac{2}{p}}} \frac{\vert \g(y)-\g(x)\vert}{\vert x-y\vert^{2}} &\leq& C \frac{\g(x)^{\frac{2}{p}-1}+\g(y)^{-\frac{2}{p}}}{\vert x-y \vert^{2}} \\
&\leq& C' \frac{1}{\vert x-y \vert^{2-\beta \max(\frac{1}{2}-\frac{1}{p}, \frac{1}{p})}}
\end{eqnarray*}

\noindent Now, let us introduce the function $x \mapsto \Theta(x)$ as follows
$$
\Theta(x) \equiv \min\left(1, \frac{1}{\vert x \vert^{2-\beta(\frac{1}{2}-\frac{1}{p}, \frac{1}{p})}   }\right)
$$
Thus, we have proved that
$$
\left \vert \frac{1}{\sqrt{\w(x)}} [\Lambda, \sqrt{\w(x)} ] f  \right\vert \leq C \frac{1}{\w(x)^{\frac{1}{p}}} \int \Theta(x-y) \w(y)^{\frac{1}{p}} \vert f(y) \vert \ dy + C \vert \mathcal{H}^{\#}f(x) \vert
$$
Since $2-\beta \max(\frac{1}{2}-\frac{1}{p}, \frac{1}{p}) > \frac{3}{2}$, then the function $\Theta$ is an integrable function on $\mathbb R$. Taking the power $p$ in both side, multiplying by $w$ and then integrating with respect to $x$ give the following
$$ 
\int  \vert \frac{1}{\sqrt{\w(x)}} [\Lambda, \sqrt{\w(x)} ] f   \vert^{p} \ \w dx \leq C \int (\Theta * G)(x) \ dx + C' \int \vert \mathcal{H}^{\#}f(x) \vert^{p} \ \w dx
$$
where we set $G(y)=\w(y)^{1/p} \vert f (y) \vert$. Therefore, since $\Theta \in L^{1}(\mathbb R)$ and $G \in L^{p}(\mathbb R)$, Young's inequality for the convolution gives
$$ 
  \left\Vert \frac{1}{\sqrt{w(x)}} [\Lambda, \sqrt{\w(x)} ] f   \right\Vert^{p}_{L^{p}(\w dx)}  \leq C'' \int \vert f (x) \vert^{p} \ \w \ dx 
  $$
where, in the second part of the inequality, we have used that the truncated Hilbert transform of $f$ is a Calder\'on-Zygmund type operator and as such is bounded on $L^{p}(\w dx)$ ( by the $L^{p}(\w dx)$ norm of $f$) since $\w \in \mathcal{A}_{p}$ for all $p \in [2, \infty)$. This concludes the proof of the second commutator estimate.
\qed

\subsection{Bounds for  $\Lambda w_{\beta}$}

We have used in the previous subsection  the bound   $\vert \partial_{x} w_{\beta}(x) \vert \leq C w_{\beta}(x)$.
A similar estimate holds for the nonlocal operator $\Lambda$ : 
\begin{lemma} \label{borne}
For all $\beta \in (0,1)$, we have $\vert \Lambda w_{\beta}(x) \vert \leq C w_{\beta}(x)$
\end{lemma}
 \noindent ${\bf{Proof \ of \ lemma \ \ref{borne}}}.$ We need to estimate the following singular integral
 $$
 \Lambda w_{\beta} (x) = \frac{P.V.}{\pi} \int \frac{w_{\beta}(x)-w_{\beta}(y)}{\vert x-y \vert^{2}} \ dy
 $$
To do so, we split the integral in three pieces   
$$
\frac{P.V.}{\pi} \int \frac{w_{\beta}(x)-w_{\beta}(y)}{\vert x-y \vert^{2}} \ dy = \frac{P.V.}{\pi}  \sum_{i=1}^{3} \int_{\Delta_{i}(x)}  \frac{w_{\beta}(x)-w_{\beta}(y)}{\vert x-y \vert^{2}} \ dy \equiv \sum_{i=1}^{3} I_{i}
$$
The domains of integrations  $\Delta_{i} (x)$ with $i=1,2,3$ are the same  ones as those introduced in the previous subsection. Using (\ref{trunc}), we get the following estimate for the integration in $\Delta_{1}(x)$
\begin{eqnarray*}
I_1&\leq&\frac{P.V.}{\pi}  \int_{\Delta_{1}(x)}  \frac{\vert w_{\beta}(x)-w_{\beta}(y) \vert }{\vert x-y \vert^{2}} \ dy \\
&\leq& \frac{P.V.}{\pi}  \int_{\Delta_{1}(x)}  \frac{\vert w_{\beta}(y)-w_{\beta}(x) + \alpha(x-y)(x-y) \partial_{x} w_{\beta}(x) \vert}{\vert x-y \vert^{2}} \ dy \\
 &\leq& C w_{\beta}(x) 
\end{eqnarray*}
For the integral over $\Delta_{2}(x)$, we have 
$$
I_2 \leq \frac{P.V.}{\pi}  \int_{\Delta_{2}(x)}  \frac{\vert w_{\beta}(x)-w_{\beta}(y) \vert }{\vert x-y \vert^{2}} \ dy \leq \frac{P.V.}{\pi}  \int_{\Delta_{2}(x)}  \frac{\vert w_{\beta}(x) \vert}{\vert x-y \vert^{2}} \ dy < C w_{\beta}(x)
$$
The last integral can be estimated as follows
\begin{eqnarray*}
\frac{P.V.}{\pi}  \int_{\Delta_{3}(x)}  \frac{\vert w_{\beta}(x)-w_{\beta}(y) \vert }{\vert x-y \vert^{2}} \ dy &\leq& C   \int_{\Delta_{3}(x)}  \frac{\vert w_{\beta}(x) \vert}{\vert x-y \vert^{2}} \ dy + C \int_{\Delta_{3}(x)}  \frac{1}{\vert x-y \vert^{2+\beta}} \ dy \\ 
&\leq& C' \w(x)
\end{eqnarray*}
This concludes the proof of the lemma.
\qed
 
\section{A priori estimates in  weighted Sobolev spaces}

In order to prove the theorems, we approximate our initial data by data which vanish at infinity, so that we may use the existence and regularity results obtained in the last section (see section \ref{H3}).  For a solution $\theta$ in $H^s$, $s=0$, $1/2$ or $1$, we have obviously $\theta\in H^s(w_\beta\, dx)$. This will allow us to estimate the norm of $\theta$ in $H^s(w_\beta\, dx)$; we shall show that those estimates do not depend on the $H^s(dx)$ norm of $\theta_0$, but only on the norm of $\theta_0$ in $H^s(w_\beta\, dx)$ and thus  we shall be able to relax the approximation. \\

 In the sequel,  we shall just write $w$ instead of $w_{\beta}$ for the   sake of readibility. 

\subsection{Estimates for the $L^{2}(\dw)$ norm}\label{l2w}

In this subsection, we consider the solution $\theta\in H^1$  associated to some   initial value $\theta_0\in H^1$ and try to estimate its $L^2(\dw)$ norm.

As usually, we multiply the transport equation by $w\theta$ and we integrate with respect to the space variable. We  obtain
\begin{eqnarray*}
\frac{1}{2} \frac{d}{dt}\left( \int \theta^{2} w \ dx \right) &=& \int \theta \partial_{t} \theta \ w \ dx \\
&=& - \int \theta \Lambda \theta \ w dx - \int \theta \mathcal{H} \theta \partial_{x} \theta w \ dx.
\end{eqnarray*}
When integrating by parts, we take into account the weight $w$ and get 
 \begin{equation*}\begin{split}  \frac{1}{2}\frac{d}{dt} \left(\int \theta^2\, dx\right)= & -\int \vert\Lambda^{1/2} \theta\vert^2\, w\,  dx-\frac{1}{2}\int \theta^2\, \Lambda\theta\, w\, dx\\ &- \int \Lambda^{1/2} \theta [\Lambda^{1/2}, w] \theta \ dx + \frac{1}{2} \int \theta^{2} \mathcal{H} \theta \ \partial_{x}w \ dx.
\end{split}\end{equation*}
Using   lemma \ref{comm}
\begin{eqnarray*}
\int  \Lambda^{1/2} \theta \ \frac{1}{w}[\Lambda^{1/2}, w] \theta\, w \ dx &\leq& \Vert \Lambda^{1/2} \theta \Vert_{L^{2}(w \ dx)} \left \Vert \frac{1}{w} [\Lambda^{1/2}, w] \theta \right\Vert_{L^{2}(wdx )} \\
&\leq& C \Vert \Lambda^{1/2} \theta \Vert_{L^{2}(w  dx)} \Vert \theta \Vert_{L^{2}(\dw)} \\
&\leq& \frac{1}{2} \int \vert \Lambda^{1/2} \theta \vert^{2} \ w \ dx \ + \ \frac{C^{2}}{2} \int \theta^{2} \ w \ dx.
\end{eqnarray*}
Moreover, we have 
\begin{eqnarray*}
\frac{1}{2} \int \theta^{2} \mathcal{H} \theta \ \partial_{x}w \ dx&\leq & C \|\theta\|_\infty  \int \vert\theta\vert \vert\mathcal{H} \theta \vert w \ dx  \\
&\leq& C'  \|\theta_0\|_\infty \|\theta\|_{L^2(\dw)}^2
\end{eqnarray*}
Thus, we find that
$$
 \frac{d}{dt}\left( \int \theta^{2} w \ dx \right) + \int  \vert \Lambda^{1/2} \theta \vert^{2} \ w \ dx  \leq C(1+\Vert \theta_{0} \Vert_{\infty}) \int \theta^{2} \ w \ dx- \int \theta^2 \Lambda \theta\, w\, dx
$$
If $\theta_0$ is nonnegative, then the maximum principle gives us that $\theta\geq 0$.
Then, using the pointwise C\'ordoba and C\'ordoba inequality \cite{CC}  (valid for   $  \theta\geq 0$)
  \begin{equation*}  \Lambda(\theta^3)\leq 3   \theta^2\Lambda\theta   \end{equation*}   and using  lemma \ref{borne}, we get
$$
\frac{1}{2} \int \theta^{2} \partial_{x} \mathcal{H} \theta \ w \ dx \leq -\frac{1}{6} \int \Lambda(\theta^{3}) \ w \ dx = -\frac{1}{6} \int \theta^{3} \Lambda w \ dx \leq C \Vert \theta_{0} \Vert_{\infty} \int \theta^{2} w \ dx
$$
Integrating in time $s\in [0,T]$ we conclude thanks to  Gronwall's lemma that we have a global control of both $\Vert \theta \Vert_{L^{\infty}([0,T], L^{2}(w dx))}$ and  $\Vert \Lambda^{1/2} \theta \Vert_{{L^{2}([0,T],L^{2}(\dw))}}$ by $\Vert \theta_{0} \Vert_{L^{2}(\dw)}$ and $\|\theta_0\|_\infty$.

\begin{remark}
If no assumption is made on the sign of $\theta_0$, we  just obtain
\begin{equation}\label{eql2}
 \frac{d}{dt}\left( \int \theta^{2} w \ dx \right) + \int  \vert \Lambda^{1/2} \theta \vert^{2} \ w \ dx  \leq C(1+\Vert \theta_{0} \Vert_{\infty}) \int \theta^{2} \ w \ dx + \|\theta_0\|_\infty \int \vert  \theta  \Lambda \theta\vert \, w\, dx,
\end{equation}
 which requires a control on $\|\Lambda\theta\|_{L^2(\dw)}$.
 \end{remark}

\subsection{Estimates for the $H^{1/2}(\dw)$ norm }
\noindent In this subsection, we consider the evolution norm of  $\theta $   in $H^{1/2}(\dw)$. We have

\begin{eqnarray*}
\frac{1}{2} \frac{d}{dt} \left( \int \vert \Lambda^{1/2}\theta\vert^2 \ wdx \right)&=& \int \partial_t \theta \Lambda^{1/2}(w \Lambda^{1/2} \theta)dx \\
&=& - \int \Lambda\theta \Lambda^{1/2} (w \Lambda^{1/2}\theta)  dx - \int  \mathcal{H}\theta \partial_{x} \theta \Lambda^{1/2}(w \Lambda^{1/2} \theta)  \ dx.
\end{eqnarray*}
Then, we get the weight $w$ outside from the differential terms
\begin{eqnarray*}
\frac{1}{2} \frac{d}{dt} \left( \int \vert \Lambda^{1/2}\theta \vert^2 \ wdx \right)&=& - \int \vert \Lambda \theta \vert^2 \ \dw  - \int \mathcal{H} \theta \partial_x \theta \Lambda \theta  \ \dw \\
&+& \int \Lambda \theta (w \Lambda^{1/2}\Lambda^{1/2} \theta- \Lambda^{1/2}(w \Lambda^{1/2}\theta)) \ dx \\
&+& \int  \left(w \Lambda^{1/2} \Lambda^{1/2} \theta - \Lambda^{1/2}(w \Lambda^{1/2} \theta) \right)\mathcal{H}\theta \partial_{x} \theta \ dx
\end{eqnarray*}
Finally, we distribute in the second term the weight $w=\gamma^2$ equally into the $\partial_x$ and the $\Lambda$ term, we obtain
\begin{eqnarray*}
\frac{1}{2} \frac{d}{dt} \left( \int \vert \Lambda^{1/2}\theta\vert^2 \ \dw \right)= &-& \int \vert \Lambda \theta \vert^2 \ \dw - \int \mathcal{H} \theta \partial_{x}(\gamma \theta) \Lambda(\gamma \theta) \ dx \\
&-& \int \mathcal{H} \theta \gamma \Lambda \theta \ (\gamma \partial_x \theta- \partial_{x}(\gamma \theta)) \ dx \\
&-& \int \partial_{x}(\gamma \theta) \mathcal{H}\theta (\gamma \Lambda \theta-\Lambda(\gamma \theta)) \ dx \\
&+& \int \Lambda \theta (w \Lambda \theta - \Lambda^{1/2}(w \Lambda^{1/2}\theta))  \ dx \\
&+&\int  (w \Lambda \theta - \Lambda^{1/2}(w \Lambda^{1/2}\theta)) \mathcal{H} \theta \partial_x \theta \ dx \\
= &-& \int \vert \Lambda \theta \vert^2 \ \dw + J_1+J_2+J_3+J_4+J_5
\end{eqnarray*}
Let us estimate $J_{1}$. Using the  $\mathcal{H}^1$-$BMO$ duality, we write  \\
$$ J_1 \leq C'_1 \Vert \mathcal{H}\theta \Vert_{BMO} \Vert \partial_{x}(\gamma \theta) \Lambda(\gamma \theta) \Vert_{\mathcal{H}^1}$$
Now, we shall use  the fact that if a function $f\in L^2$ then the function $g=f \mathcal{H}f$ belongs to the Hardy space $\mathcal{H}^{1}$ : indeed, we have 

\begin{equation} \label{magic}
2 \mathcal{H} (f \mathcal{H}f)(x)=(\mathcal{H}f(x))^2  - f(x)^2
\end{equation}  so that $f\mathcal{H}f$ belongs to $\mathcal{H}^1$ and we have
$$
\Vert f \mathcal{H}f \Vert_{\mathcal{H}^1}=\Vert f \mathcal{H}f \Vert_1+\Vert  \mathcal{H}\big(f \mathcal{H}f \big)\Vert_1 \leq C \Vert f \Vert^{2}_{L^2}
$$
From    formula (\ref{magic}), we get the following estimate
$$
J_{1} \lesssim \Vert \theta_{0} \Vert_{\infty} \Vert \partial_{x}({\gamma \theta}) \Vert^{2}_{L^{2}} \lesssim   \Vert \theta_{0} \Vert_{\infty} \left( \Vert  \theta\Vert^{2}_{L^{2}(wdx)} + \Vert \Lambda{ \theta} \Vert^{2}_{L^{2}(wdx)}\right).
$$
To estimate $J_{2}$, we use the fact that $\vert \partial_x \gamma \vert < C'_{2} \gamma$ and that $w \in \mathcal{A}_{4}$,  we obtain 
\begin{eqnarray*}
J_2 = \int \mathcal{H}\theta \  \gamma \Lambda \theta \ \theta\partial_x \gamma  \ dx &\lesssim&  \int  \left\vert w^{1/2}\mathcal{H}\theta \  w^{1/2} \Lambda \theta \ \theta\right\vert  \ dx\\
 &\lesssim& C \Vert \mathcal{H\theta} \Vert_{L^{2}(wdx)} \Vert \Lambda \theta \Vert_{L^{2}(w dx)} \Vert \theta \Vert_{L^{\infty}}.
 \end{eqnarray*}
Therefore,
$$
J_2 \lesssim \Vert \theta_{0} \Vert_{\infty} \Vert \theta \Vert_{L^{2}(wdx)} \Vert \Lambda \theta \Vert_{L^{2}(w dx)}.
$$

\noindent In order to estimate $J_{3}$, we take $p_1$ and $q_1$ with $2<p_1<\infty$ and $\frac{1}{p_1} + \frac{1}{q_1}=\frac{1}{2}$ and using  lemma \ref{comm} we obtain
\begin{eqnarray*}
J_{3} &\leq& \Vert \partial_{x}(\gamma \theta) \Vert_{2} \Vert \mathcal{H}\theta (\gamma \Lambda \theta-\Lambda(\gamma \theta)) \Vert_2 \\
&\lesssim& \Vert \partial_{x}(\gamma \theta) \Vert_2 \Vert \mathcal{H} \theta \Vert_{L^{q_1}( \dw)} \left\Vert \frac{1}{\gamma} (\gamma \Lambda \theta- \Lambda(\gamma \theta) \right\Vert_{L^{p_1}(\dw)} \\
&\lesssim&  \Vert \partial_{x}(\gamma \theta) \Vert_2 \Vert \theta \Vert_{L^{q_1}(\dw)} \Vert \theta \Vert_{L^{p_1}(\dw)}
\end{eqnarray*}
Then, using
$$
\Vert \theta \Vert_{L^{r}(\dw)} \leq C \Vert \theta \Vert^{1-\frac{2}{r}}_{\infty} \Vert \theta \Vert^{\frac{2}{r}}_{L^{2}(\dw)}
$$
with $r=p_1$ and $r=q_1$, we find,
$$
J_3 \lesssim \Vert \theta_{0} \Vert_{\infty} \Vert \theta\Vert_{L^{2}(w dx)} \left(\Vert \theta\Vert_{L^{2}(w dx)} + \Vert \Lambda \theta \Vert_{L^{2}(w  dx)}\right)
$$
The estimation of $J_{4}$ is easy, it suffices to use lemma \ref{comm}
$$
J_4 \leq \Vert \Lambda \theta \Vert_{L^{2}(w dx)} \left\Vert \frac{1}{w}[\Lambda^{1/2},w] \Lambda^{1/2} \theta \right\Vert_{L^{2}(wdx)} \leq C_4 \Vert \Lambda \theta \Vert_{L^{2}(wdx)} \Vert \Lambda^{1/2} \theta \Vert_{L^{2}(w dx)}
$$
It remains to estimate $J_{5}$. We first write $1=w^{1/2}w^{-1/2}$ and use Cauchy-Schwarz inequality. Then, we use H\"older's inequality with $\frac{1}{p} + \frac{1}{q}= \frac{1}{2}$. We take $p$ and $q$ with $2<p<4$, and we assume $p$ to be close enough to 2 to grant that $\frac{3}{2}-\beta(1-\frac{1}{p}) > 1$ so that we may apply lemma \ref{comm}. Therefore, we get
\begin{eqnarray*}
J_5 &\leq& \Vert w^{1/2}\partial_x \theta \Vert_{L^{2}} \left\Vert w^{-1/2}\mathcal{H}\theta [\Lambda^{1/2},w] \Lambda^{1/2} \theta \right\Vert_{L^{2}} \\
&\lesssim&   \Vert \partial_x \theta \Vert_{L^{2}(wdx)}  \Vert \mathcal{H}\theta \Vert_{L^{q}(wdx)}  \left\Vert  \frac{1}{w}[\Lambda^{1/2},w] \Lambda^{1/2} \theta \right\Vert_{L^{p}(wdx)} \\
&\lesssim&   \Vert \partial_x \theta \Vert_{L^{2}(wdx)}  \Vert \theta \Vert_{L^{q}(wdx)} \Vert \Lambda^{1/2} \theta \Vert_{L^{p}(wdx)} 
\end{eqnarray*}
Moreover, using the following weighted Gagliardo-Nirenberg inequality (see inequality \ref{GN})
$$
 \Vert \Lambda^{1/2} \theta \Vert_{L^{4}(w dx)} \lesssim \Vert \theta \Vert^{1/2}_{\infty} \Vert \Lambda \theta \Vert^{1/2}_{L^{2}(wdx)},
$$
we get
\begin{eqnarray*}
\Vert \Lambda^{1/2} \theta \Vert_{L^{p}(wdx)} &\leq& \Vert \Lambda^{1/2} \theta \Vert_{L^{4}(w dx)}^{2-\frac{4}{p}} \Vert \Lambda^{1/2} \theta \Vert^{\frac{4}{p}-1}_{L^{2}(wdx)} \\
&\lesssim&  \Vert \theta_0 \Vert^{1-\frac{2}{p}}_{\infty} \Vert \Lambda \theta \Vert^{1-\frac{2}{p}}_{L^{2}(w dx)} \Vert \Lambda^{1/2} \theta \Vert^{\frac{4}{p}-1}_{L^{2}(wdx)}.
\end{eqnarray*}
Then, since
$$
\Vert \theta \Vert_{L^{q}(wdx)} \leq  \Vert \theta_0 \Vert^{\frac{2}{p}}_{\infty} \Vert  \theta \Vert^{1-\frac{2}{p}}_{L^{2}(w dx)},
$$
we get
$$
J_5 \lesssim  \Vert \theta_0 \Vert_{L^{\infty}} \Vert \theta \Vert^{1-\frac{2}{p}}_{L^{2}(w dx)} \Vert \Lambda \theta \Vert^{2-\frac{2}{p}}_{L^{2}(w dx)} \Vert \Lambda^{1/2} \theta \Vert^{\frac{4}{p}-1}_{L^{2}(wdx)}.
$$
Using the fact that
$$
\Vert \Lambda^{1/2}\theta \Vert^{\frac{4}{p}-1}_{L^{2}(wdx)} \leq  \Vert \theta \Vert^{\frac{2}{p}-\frac{1}{2}}_{L^{2}(w dx)}   \Vert \Lambda \theta \Vert^{\frac{2}{p}-\frac{1}{2}}_{L^{2}(w dx)}, 
$$
we obtain
$$
J_{5} \leq C_5 \Vert \theta_0 \Vert_{L^{\infty}} \Vert \theta \Vert^{1/2}_{L^{2}(w dx)}  \Vert \Lambda \theta \Vert^{3/2}_{L^{2}(w dx)}.
$$
Using Young's inequality, we finally   find that there exists   constants $C_{6}>0$ and $C_7>0$ (where $C_7$  depends  on $\Vert \theta_{0} \Vert_{\infty}$), such that
\begin{equation}\label{eqsob2}\begin{split}
\frac{d}{dt} \int \vert \Lambda^{1/2} \theta \vert^2 \ \dw  \leq& (C_{6} \Vert \theta_{0} \Vert_{\infty}-1) \int \vert \Lambda \theta \vert^{2}  \ \dw \\
 &+ C_7 \left(\int \theta^2 \ \dw + \int \vert\Lambda^{1/2}\theta\vert^2 \ \dw\right).
 \end{split}
\end{equation}
Combining (\ref{eql2}) and (\ref{eqsob2}), we finally obtain  \begin{equation}\label{eqsob3}\begin{split}
\frac{d}{dt}\left( \int \vert \theta\vert^2+ \vert \Lambda^{1/2} \theta \vert^2 \ \dw \right)  \leq& (C_{8} \Vert \theta_{0} \Vert_{\infty}-1) \int \vert \Lambda \theta \vert^{2}  \ \dw \\
 &+ C_9 \left(\int \theta^2 \ \dw + \int \vert\Lambda^{1/2}\theta\vert^2 \ \dw\right)\end{split}
\end{equation}  By Gronwall's lemma, we conclude that we have a control of $\Vert  \theta \Vert_{L^{\infty} L^{2}(wdx)}$, of $\Vert \Lambda^{1/2} \theta \Vert_{L^{\infty} L^{2}(wdx)}$ and of $\Vert  \Lambda\theta \Vert_{L^2L^2(wdx)}$ by $\Vert \theta_0\Vert_\infty$,  $\Vert \theta_0 \Vert_{L^2(wdx)}$ and $\Vert \Lambda^{1/2} \theta_0 \Vert_{L^{2}(wdx)}$ (if $\Vert \theta_{0} \Vert_{\infty} < \frac{1}{C_{8}}$, where $C_{8}>0$ is a constant depending only on $\beta$).

\subsection{Estimates for the $H^{1}(\dw)$ norm}  

In this subsection, we  estimate the norm of $\theta$ in $H^{1}(\dw)$. \\

In order to study the evolution of the $H^{1}( \dw)$ norm  of $\theta$, we shall study the evolution of the semi-norm $\Vert \partial_x \theta \Vert_{L^{2}(w dx)}$ instead of $\Vert \Lambda \theta \Vert_{L^{2}(w dx)}$ since they are equivalent (see Remark \ref{h1}). Therefore, we write
\begin{eqnarray*}
\frac{1}{2} \frac{d}{dt} \left( \int \vert \partial_{x} \theta \vert^{2} \ wdx \right)&=&-\int \partial_{t} \theta \ \partial_{x}(w \partial_{x}\theta) \ dx \\
&=&\int (\partial_{x} \theta)^{2} \  \mathcal{H} \theta \  \partial_{x} w  \ dx +\int \partial_{x}\theta \ \Delta \theta \ \mathcal{H} \theta  \,   w  \ dx \\
 && \ +  \int \Lambda \theta \  \partial_{x}\theta \ \partial_{x}w \ dx + \int \Lambda \theta \Delta \theta \ w \ dx   \\
\end{eqnarray*}
The last term  which come from the  linear part of the equation can be rewritten as 
$$
\int \Lambda \theta \Delta \theta \ w \ dx = -\int \Lambda \theta \Lambda^{2} \theta \ w \ dx = - \int \Lambda^{3/2} \theta [\Lambda^{1/2}, w] \Lambda \theta - \int \vert \Lambda^{3/2} \theta \vert^{2} \ w \ dx
$$
Moreover, an integration by parts gives 
$$
\frac{1}{2}\int (\partial_{x} \theta)^{2}  \ \mathcal{H} \theta \  \partial_{x} w  \ dx = - \int \partial_{x}\theta  \ \Delta \theta \ \mathcal{H} \theta \ w \ dx - \frac{1}{2}\int (\partial_{x} \theta)^{2} \ \Lambda \theta \ w \ dx
$$
So that, we get
 \begin{eqnarray*}
\frac{1}{2} \frac{d}{dt} \left( \int \vert \partial_{x} \theta \vert^{2} \ wdx \right)&=& -\int \vert \Lambda^{3/2} \theta \vert^{2} \ w \ dx - \int \Lambda^{3/2} \theta [\Lambda^{1/2}, w] \Lambda \theta  -\frac{1}{2} \int (\partial_{x} \theta)^{2} \Lambda \theta \ w \ dx \\
 &+& \frac{1}{2}\int (\partial_{x} \theta)^{2}  \ \mathcal{H} \theta \  \partial_{x} w  \ dx + \int \partial_{x}\theta\Lambda \theta \ \partial_{x}w \ dx \\
&=& - \int \vert \Lambda^{3/2} \theta \vert^{2} \ w \ dx  + J_1 + J_2 + J_3 + J_4 
\end{eqnarray*}
To estimate $J_1$ we write
$$
J_{1} = -\int  w(x) \Lambda^{3/2} \theta\  \frac{1}{w(x)}[\Lambda^{1/2}, w] \Lambda \theta \ dx \leq \Vert  \Lambda^{3/2} \theta \Vert_{L^{2}(wdx)} \Vert \frac{1}{w(x)}[\Lambda^{1/2}, w] \Lambda \theta \Vert_{L^{2}(w dx)}
$$
Therefore, using the second part of $\ref{comm}$, we conclude that
$$
J_1 \leq C_1 \Vert \Lambda^{3/2} \theta \Vert_{L^{2}(wdx)}  \Vert \Lambda \theta \Vert_{L^{2}(w dx)}
$$
For $J_{2}$, using Holder's inequality together with the fact that $w_{\beta} \in \mathcal{A}_{3}$ allows us to get
$$
J_{2}= -\frac{1}{2} \int (\partial_{x} \theta)^{2}  \Lambda\theta \ w \ dx = -\frac{1}{2} \int w^{\frac{1}{3}}\partial_{x} \theta \ w^{\frac{1}{3}}\partial_{x} \theta \ w^{\frac{1}{3}} \mathcal{H}\partial_{x} \theta \ dx \leq C \Vert \partial_{x}\theta \Vert^{3}_{L^{3}(wdx)}
$$
Then, using the following weighted Gagliardo-Nirenberg inequality
$$
\Vert \partial_{x} \theta \Vert_{L^{3}(wdx)}\leq C_2 \Vert \theta \Vert^{1/3}_{\infty} \Vert \Lambda^{3/2} \theta \Vert^{2/3}_{L^{2}(wdx)}
$$
we get
$$
J_{2} \leq C_2  \Vert \theta \Vert_{\infty} \Vert \Lambda^{3/2} \theta \Vert^{2}_{L^{2}(wdx)}
$$
The estimation of $J_{3}$ and $J_{4}$ are quite similar to the estimation of $J_2$. Indeed, we have $$
 J_{3} \leq C'_3    \int  (\partial_{x} \theta)^{2}  \ \vert \mathcal{H} \theta \vert  \   w \  dx 
  \leq  C_3 \Vert \partial_{x}\theta \Vert^{2}_{L^{3}(wdx)} \Vert \theta \Vert_{L^{3}(wdx)} 
$$
Then, using the interpolation inequality
$$
\Vert \theta \Vert_{L^{3}(wdx)} \leq \Vert \theta \Vert^{1/3}_{\infty} \Vert \theta \Vert^{2/3}_{{L^{2}(w dx)}},
$$
together with the  Gagliardo-Nirenberg inequality previously recalled, we get
$$
J_{3} \leq  C_3  \Vert \theta \Vert_{\infty} \Vert \Lambda^{3/2} \theta \Vert^{4/3}_{L^{2}(wdx)} \Vert \theta \Vert^{2/3}_{{L^{2}(w dx)}}
$$
For $J_4$, we write
$$
J_{4} \leq  C'_4 \int w^{\frac{1}{2}}\vert \partial_{x}\theta \vert \ w^{\frac{1}{2}} \vert \mathcal{H}\partial_{x} \theta \vert \  \ dx \leq C_4 \Vert \partial_{x}\theta \Vert^{2}_{L^{2}(wdx)}
$$
Therefore, by the maximum principle for the $L^{\infty}$ norm and Young's inequality, we get 
\begin{eqnarray*}
\frac{1}{2} \frac{d}{dt} \left( \int \vert \partial_{x} \theta \vert^{2} \ wdx \right) &\leq& (C_2\Vert \theta_0 \Vert_{\infty}-1)  \int \vert \Lambda^{3/2} \theta \vert^{2} \ w \ dx \\ 
&+& C_1 \Vert \Lambda^{3/2} \theta \Vert^{2}_{L^{2}(wdx)}  \Vert \theta \Vert_{L^{2}(w dx)} + C_4 \Vert \partial_{x}\theta \Vert^{2}_{L^{2}(wdx)} \  \\
&+& C_3  \Vert \theta \Vert_{\infty} \Vert \Lambda^{3/2} \theta \Vert^{4/3}_{L^{2}(wdx)} \Vert \theta \Vert^{2/3}_{{L^{2}(w dx)}}\\
&\leq& ({C'_2}  \Vert \theta_0 \Vert_{\infty}-1) \int \vert \Lambda^{3/2} \theta \vert^{2} \ \dw  \\
 && \ + C_5 \left( \Vert \theta \Vert^{2}_{L^{2}(wdx)} +  \Vert \partial_{x}\theta \Vert^{2}_{L^{2}(wdx)}\right),
\end{eqnarray*} where the constant $C_5$ depends on $\Vert\theta_0\Vert_\infty$. Then, integrating in time $s\in[0,T]$ gives
\begin{eqnarray}
\Vert \theta(T,.) \Vert^{2}_{H^{1}(wdx)} &\leq& (C'_{2}\Vert \theta_0 \Vert_{\infty}-1) \int_{0}^{T}  \Vert \Lambda^{3/2} \theta \Vert^{2}_{L^{2}(wdx)} \ ds \nonumber \\
&& + \ C_5  \int_{0}^{T}  \Vert \theta(s,.) \Vert^{2}_{H^{1}(wdx)} \ ds
\end{eqnarray}

Therefore, Gr\"onwall's lemma allows us to conclude that we have a global control of  $\Vert \partial_{x}\theta\Vert_{L^{\infty}L^{2}(w dx)}$ and  $\Vert \Lambda^{3/2} \theta \Vert_{L^{2}L^{2}(wdx)}$ by $ \Vert \theta_0 \Vert_{\infty}$ and $\Vert \theta_{0} \Vert_{H^{1}(wdx)}$, provided that $\Vert \theta_0 \Vert_{\infty} < \frac{1}{{C'_2}}$. Note that $C'_2>0$ is a constant that depends only on $\beta$.

\section{Proof of the theorems}
\subsection{The truncated initial data}

We shall approximate $\theta_0$ by $\theta_{0,R}=\theta_0(x)\psi(\frac{x}{R})$, where $\psi$ satisfies the following assumptions :
\begin{itemize}
\item $\psi\in \mathcal{D}(\mathbb{R})$
\item $0\leq \psi\leq 1$
\item $\psi(x)=1$ for $x\in [-1,1]$ and $=0$ for $\vert x\vert\geq 2$
\end{itemize}

This approximation neither alters the non-negativity of the data, nor increases its $L^{\infty}$ norm.
We have obviously the strong convergence, when $R\rightarrow +\infty$,  of $\theta_{0,R}$ to $\theta_0$ in $H^s(w\, dx)$ if $\theta_0\in H^s(w\, dx)$ and $s=0$ or $s=1$. The only difficult case is $s=1/2$. This could be dealt with through an interpolation argument. But we shall give a direct proof that $$\lim_{R\rightarrow +\infty} \|\Lambda^{1/2}(\theta_0-\theta_{0,R})\|_{L^2(w\, dx)}=0.$$ As we have the strong convergence of $\psi_R \Lambda^{1/2}\theta_0$ to $\Lambda^{1/2}\theta_0$ in $L^2(w\, dx)$, we must  estimate the norm of the commutator $[\Lambda^{1/2},\psi_R]\theta_0$ in $L^2(w\, dx)$, where we write $\psi_R(x)=\psi(\frac{x}{R})$.
We just write
$$
\left \vert [\Lambda^{1/2}, \psi_R] \theta_{0} \right \vert \leq C \int \frac{\vert \psi_R(x)-\psi_{R}(y) \vert}{\vert x-y \vert^{3/2}}  \vert \theta_{0}(y) \vert \ \ dy
$$
with
$$
\frac{\vert \psi_R(x)-\psi_{R}(y) \vert}{\vert x-y \vert^{3/2}}\leq   \min\left( \frac{\Vert \partial_{x} \psi \Vert_{\infty}}{R \vert x-y \vert^{1/2}}, \frac{2\Vert \psi \Vert_{\infty}}{\vert x-y \vert^{3/2}} \right)=  \frac{1}{R^{3/2}} K(\frac{x-y}{R})
$$
where the kernel $K$ is   integrable,  nonnegative and radially decreasing; thus, from inequality (\ref{St}), we find 
that
 $$\left\vert [\Lambda^{1/2}, \psi_R] \theta_{0} \right\vert\leq \|K\|_1  R^{-1/2}\mathcal{M}\theta_0$$ which gives
 $$ \Vert  [\Lambda^{1/2}, \psi_R] \theta_{0}\Vert_{L^2(w\, dx)}\leq C R^{-1/2} \|\theta_0\|_{L^2(w\, dx)}.$$

\subsection{Proof of theorem \ref{th22}}
 
We consider the sequence $ \theta_{0,N}$, $N\in\mathbb{N}$ and $N\geq 1$.  We have the convergence of $\theta_{0,N}$  to $\theta_0$ in $H^{1/2}(w\, dx)$. Moreover, if $\|\theta_0\|_\infty$ is small enough we know that we have a solution $\theta_N$ of our transport equation $\mathcal{T}$ with initial value $\theta_{0,N}$.    Using the {\it{a  priori}} estimates of the previous section, we get  (uniformly with respect to $N$) that the sequence $\theta_N$ is bounded in the space $L^{\infty}([0,T], H^{1/2}(\dw))$ and $L^{2}([0,T], H^{1}(\dw))$  for every $T \in  (0,\infty)$. Now, let $\psi(x,t) \in \mathcal{D}((0,\infty] \times \mathbb R)$, then $\psi\theta_N$ is bounded in $L^{2}([0,T], H^{1})$. Moreover, we have 
$$\partial_{t}({\psi}\theta_N)=\theta_N\partial_{t}\psi  + \psi \partial_{t}\theta_N=(I)+(II)$$
Obviously, $(I)$ is bounded in $L^{2}([0,T], L^{2})$. For $(II)$, we write
$$
 \psi \partial_{t}\theta_N= -\psi\partial_{x}{\theta_N} \mathcal{H}\theta_N - \psi\Lambda\theta_N= 
 -\psi \partial_{x}(\theta_N \mathcal{H}\theta_N) + \psi\theta_N \Lambda\theta_N - \psi\Lambda\theta_N
$$
Since $\theta_N$ is bounded in $L^{2}([0,T], L^{2}(w\, dx))$ then by the continuity of the Hilbert transform on $L^{2}$, the sequence $\mathcal{H}\theta_N$ is bounded in $L^{2}([0,T], L^{2} (w\, dx))$ therefore, since $\theta_N$ is bounded in $L^{\infty}([0,T], L^{\infty})$, we get  that $\psi \partial_{x}(\theta_N \mathcal{H}\theta_N)$ (=$\partial_{x}(\psi \theta_N \mathcal{H}\theta_N)  -(\partial_x\psi) \theta_N \mathcal{H}\theta_N$ )  is bounded in  $L^{2}([0,T], H^{-1})$. Therefore,   since $\psi (1-\theta_N) \Lambda\theta_N$ is bounded in $L^{2}([0,T], L^2))$ we conclude that $\partial_{t}({\psi}\theta_N)$ is bounded in $L^{2}([0,T], H^{-1})$. By Rellich compactness theorem \cite{PGLR}, there exists a subsequence $\theta_{N_{k}}$ and a function $\theta$ such that 
$$
\theta_{N_{k}} \xrightarrow[N_k\to+\infty]{} \theta \ \  \text{strongly} \ \ \text{in} \  L^{2}_{loc} ((0,\infty) \times \mathbb R),
 $$
Futhermore, since the sequence $ \theta_{N_{k}}$ is bounded in spaces whose dual space are separable Banach spaces, we get the two following *-weak convergences, for all $T<\infty$
$$
\theta_{N_{k}} \xrightarrow[N_k\to+\infty]{} \theta \ \ \  \text{*-weakly} \ \ \text{in} \  L^{\infty}([0,T], H^{1/2}(\dw)),
 $$
 and,
 $$
\theta_{N_{k}} \xrightarrow[N_k\to+\infty]{} \theta \ \ \ \text{*-weakly} \ \ \text{in} \  L^{2}([0,T], H^{1}(\dw)),
 $$
It remains to check that $\theta$ is a solution of the transport equation $\mathcal{T}$. Let $\Phi$ be a compactly supported smooth function, we need to prove the equality
$$
\int \int_{t>0} \theta \ \partial_{t} \Phi \ dx \ dt = \int \int_{t>0} \Phi \left( \mathcal{H} \theta \partial_{x} \theta + \Lambda \theta \right) \ dx \ dt - \int \Phi(0,x) \theta_{0}(x) \ dx.
$$
To prove this equality, it suffices to prove that we can pass to the weak limit in the following equality
$$
\int \int_{t>0}  \theta_{N_{k}} \ \partial_{t} \Psi \ dx \ dt = \int \int_{t>0} \Psi \left( \mathcal{H} \theta_{N_{k}}  \partial_{x}  \theta_{N_{k}} + \Lambda  \theta_{N_{k}} \right) \ dx \ dt - \int \Psi(0,x) \theta_{N_{k},0}(x) \ dx.
$$
The *-weak convergence of $ \theta_{N_{k}}$ toward $\theta$  in $L^{\infty}((0,T),L^{2}))$  implies  the convergence in  $\mathcal{D}'([0,T] \times \mathbb R)$ and therefore
$$
 \partial_t  \theta_{N_{k}}  \xrightarrow[N_k\to+\infty]{}  \partial_{t}\theta \ \ \text{in} \ \mathcal{D}'([0,T] \times \mathbb R).
$$
Moreover, since $\Lambda  \theta_{N_{k}}$ is a (uniformly) bounded sequence on $L^{2}([0,\infty] \times \mathbb R)$ therefore we also have convergence in the sense of distribution
$$
\Lambda \theta_{N_{k}} \xrightarrow[n_k\to+\infty]{} \Lambda\theta \ \ \text{in} \ \mathcal{D}'([0,T] \times \mathbb R).
$$
It remains to treat the nonlinear term, we rewrite it as
$$
\int \int_{t>0} \Psi \mathcal{H} \theta_{N_{k}}  \ \partial_{x}  \theta_{N_{k}} \ dx \ dt = -  \int \int_{t>0}  \theta_{N_{k}} \mathcal{H} \theta_{N_{k}}  \partial_x \Psi - \int \int_{t>0} \Psi  \theta_{N_{k}} \partial_{x}  \mathcal{H} \theta_{N_{k}} \ dt \ dx.
$$
Using the strong convergence of $ \theta_{N_{k}}$   on $L^{2}_{loc} ((0,\infty) \times \mathbb R)$   and the *-weak convergence of $\mathcal{H} \theta_{N_{k}}$ in  $L^{2}([0,T], L^2)$, we conclude that the products $ \theta_{N_{k}} \mathcal{H} \theta_{N_{k}}$ converge weakly in $L^{1}_{loc} ((0,\infty) \times \mathbb R)$ toward $\theta \mathcal{H}\theta$. For the second term, we also use the strong  $L^{2}_{loc} ((0,\infty) \times \mathbb R)$ convergence of $ \theta_{N_{k}}$  and the weak convergence of  $\partial_{x} \mathcal{H} \theta$ on $L^{2} ((0,\infty) \times \mathbb R)$, we conclude that the product converges in $L^{1}_{loc} ((0,\infty) \times \mathbb R)$. 
\qed \\

\subsection{Proof of theorem \ref{th23}}
The proof of Theorem \ref{th23} is   similar to the proof of Theorem \ref{th22}, using {\it{a priori}} estimates on the $H^1(w\, dx)$ norm instead  of the $H^{1/2}(w\, dx)$ norm.

\subsection{The case of data in $L^2(dx)$ or $L^2(\dw)$} \label{l2}
 When $\theta_0\in L^2\cap L^\infty$ and is non-negative, we have a priori estimates on the $L^2$ norm of $\theta$ that involves only $\|\theta_0\|_2$ and $\|\theta_0\|_\infty$, but this is not sufficient to grant existence of the solution $\theta$, as we have not enough regularity to control the nonlinear term $\mathcal{H} \theta
\partial_x\theta$.

Indeed, we have a control of $\mathcal{H}\theta$ in $L^2 H^{1/2}$ and of $\partial_x\theta $ in $L^2 H^{-1/2}$. But to pass to the limit in our use of the Rellich  theorem, we should have (local) strong convergence of $\theta_{\eta_k}$ to $\theta$ in $L^2 H^{1/2}$ while we may establish only the *-weak convergence.  This can be seen as follows : if  $\theta_n$ is a bounded sequence in $L^2 H^{1/2}$ that converge locally strongly in $L^2 L^2$ to a limit $\theta$ and if  $\mathcal{H} \theta_n
\partial_x\theta_n$ converges in $\mathcal{D}'$, we write
\begin{equation*}\begin{split}
\mathcal{H} \theta_n
\partial_x\theta_n=& \partial_x(\theta_n \mathcal{H}\theta_n)-\theta_n \partial_x\mathcal{H}\theta_n\\=&  \partial_x(\theta_n \mathcal{H}\theta_n)+ \theta_n \Lambda\theta_n  \\=& \partial_x(\theta_n \mathcal{H} \theta_n)+\frac{1}{2}\Lambda(\theta_n^2)+ C \int \frac{ (\theta_n(t,x)-\theta_n(t,y))^2}{\vert x-y\vert^2}\, dy.
\end{split}\end{equation*}
While we have the convergence in $\mathcal{D}'$ of  $\partial_x(\theta_n \mathcal{H} \theta_n)+\frac{1}{2}\Lambda(\theta_n^2)$ to $\partial_x(\theta  \mathcal{H}\theta)+\frac{1}{2}\Lambda(\theta^2)$, we can only write
$$ \lim_{n\rightarrow +\infty}  \int \frac{ (\theta_n(t,x)-\theta_n(t,y))^2}{\vert x-y\vert^2}\, dy= \int \frac{ (\theta(t,x)-\theta(t,y))^2}{\vert x-y\vert^2}\, dy+\mu,$$ where $\mu$ is a non-negative measure.

\section{The construction of regular enough solutions revisited \label{H3}}

 The global existence results of   C\'ordoba, C\'ordoba and Fontelos in \cite{CCF}  and of Dong  in \cite{Dong} correspond to  Theorems \ref{th22} to \ref{th23} in the case $\beta=0$  : they are mainly based on the maximum principle (if $\theta_0$ is bounded, then $\theta$ remains bounded and if $\theta_0$ is non-negative, $\theta$ remains non-negative) along with the use of some useful identities or inequalities involving the nonlocal operators $\Lambda$ and $\mathcal{H}$.  We do not know whether our solutions become smooth (this is known in the case $\beta=0$ for Theorem \ref{th22},  this is proved by Kiselev \cite{Ki}). Another interesting question is whether we have eventual regularity in the sense of \cite{S} for our solutions. \\

In this section, for conveniency, we sketch a complete proof of Theorems \ref{th22} and \ref{th23} in the case $\beta=0$, under a smallness assumption on $\|\theta_0\|_\infty$ (although this latter case is treated in \cite{CCF}, we shall give a slightly different proof for the {\it{a priori}} estimates). Before starting the {\it{a priori}} estimates, one has to deal with the existence issue, namely, proving the existence of at least one solution. This step is rather important for this model since for instance one can derive a nice energy estimate for the $L^2$ (resp weighted $L^2$) norm (see \cite{CCF}, resp see section \ref{l2w}) whereas the existence of such a solution  is not clear in both cases (see section \ref{l2}). Since we aim at proving global existence results and not only {\it{a priori}} estimates, we need to give a proof of the existence of regular enough solutions. This is done in six steps and is based on classical arguments.   \\

\noindent{\bf First step :  regularizations of the  equation and of the data  }\\

We use a nonnegative smooth  compactly supported function $\varphi$ (with $\int\varphi(x) \, dx=1$) and for positive  parameters  $\epsilon$, $\eta$ we consider  the parabolic approximation of  equation $(\mathcal{T}_1)$ :
    \begin{equation} 
\ (\mathcal{T}_{1}^{\epsilon,\eta}) \ : \\\left\{
\aligned
&\partial_{t}\theta+\theta_x \mathcal{H}\theta+ \nu \Lambda \theta = \epsilon\Delta \theta \hspace{2cm} 
\\ \nonumber
& \theta(0,x)=\theta_{0}*\varphi_\eta(x). \text{ (with } \varphi_\eta(x)\equiv\frac{1}{\eta} \varphi(\frac{x}{\eta})\text{ )}
\endaligned
\right.
\end{equation}
Recall that $\Delta\theta=\partial^2_x\theta$, then we can rewrite the problem into an integral form as follows
$$ \theta=e^{\epsilon t\Delta }( \theta_0*\varphi_\eta)-\int_0^t e^{\epsilon (t-s)\Delta } (\theta_x \mathcal{H}\theta+ \nu \Lambda \theta)\, ds.$$

We  may  solve this equation in   $\mathcal{C}([0,T_{\epsilon,\eta}], H^3)\cap L^2((0,T_{\epsilon,\eta}),  H^4)$, for some small enough time $T_{\epsilon,\eta}$. Indeed,  we have, for $T>0$ and  for a constant $C_\epsilon$ independent of $T$,  for all $\gamma_0\in H^3$,    $u,v\in  \mathcal{C}([0,T], H^3)\cap L^2((0,T),\dot H^4)$ and $w\in L^2([0,T],H^2)$~:
 \begin{itemize}
\item   $\displaystyle \sup_{0<t<T} \| e^{\epsilon t\Delta} \gamma_0\|_{H^3}\leq    \|\gamma_0\|_{H^3}$ and $ \|\Delta   e^{\epsilon t\Delta} \theta_0\|_{L^2((0,T), L^2)}\leq C_\epsilon \|\theta_0\|_{H^3}$
\item $\displaystyle\int_0^t e^{\epsilon(t-s)\Delta} w\, ds\in \mathcal{C}([0,T],H^3)\cap L^2([0,T], H^4)$
 \item $\displaystyle\sup_{0<t<T} \left\| \int_0^t e^{\epsilon(t-s)\Delta } w \, ds\right\|_{2}\leq C_\epsilon T^{1/2}  \| w\|_{L^2 H^2}$ 
 \item $\displaystyle\sup_{0<t<T} \left\| \partial_x^3\int_0^t e^{\epsilon(t-s)\Delta } w\, ds\right\|_{2}\leq C_\epsilon  \| w\|_{L^2 H^2}$ 
 \item  $ \left\|\Delta\int_0^t e^{\epsilon(t-s)\Delta} w\, ds\right\|_{L^2((0,T), L^2)}\leq C_\epsilon    \| w\|_{L^2  H^2}$
\item  $\displaystyle  \| \Lambda  u\|_{L^2 H^2}\leq C T^{1/2}  \| u\|_{L^\infty H^3}$ 
\item    $\displaystyle\| u_x \mathcal{H} v \|_{L^2 H^2}\leq C    T^{1/2}\|u\|_{L^\infty H^3}  \|v\|_{L^\infty  H^3}$ 
\end{itemize}
Thus, using Picard's iterative scheme, we find a solution  $$ \theta=e^{\epsilon t\Delta } \gamma_0-\int_0^t e^{\epsilon (t-s)\Delta } (\theta_x \mathcal{H}\theta+ \nu \Lambda \theta)\, ds$$ on an interval $[0,T_{\epsilon,\eta}]$, where $T_{\epsilon,\eta}$ depends only on $\epsilon$ and  $\| \gamma_0\|_2$.  If $\|\theta\|_{H^3}$ remains bounded, we may bootstrap the estimates to get an extension to a larger interval. Thus, if $T^*_{\epsilon,\eta}$ is the maximal existence time, we must have  
$$T^*_{\epsilon,\eta}<+\infty \Rightarrow   \sup_{0<t<T^*_{\epsilon,\eta}} \|\theta_{\epsilon,\eta}(t,.)\|_{H^3}=+\infty.$$

The strategy is then to have a criterion on $\theta_0$ to ensure that $T^*_{\epsilon,\eta}=+\infty$ for every $\epsilon>0$ and to get uniform controls on the solutions $\theta_{\epsilon,\eta}$ to allow to get a limit when $\epsilon$  and $\eta$ go  to $0$.

\noindent{\bf Second step :  applying the maximum principle }\\

This point is classical.  
If $\theta$ is the solution of  equation $ (\mathcal{T}^{\epsilon,\eta}) $, we define $M(t)=\displaystyle\sup_{x\in\mathbb{R}^3}  \theta(t,x) $ and $m(t)=\displaystyle\inf_{x\in\mathbb{R}^3}  \theta(t,x)  $. 
  For $t=t_0$, if $M(t_0)>0$ then the  supremum is attained at some point $x_0$, and we have $\partial_t \theta(t_0,x_0)\leq 0$, since $\Lambda \theta(t_0,x_0)\geq 0$, $\Delta \theta(t_0,x_0)\leq 0$ and $\partial_x\theta(t_0,x_0)=0$ (recall that $\theta(t_0,.)$ is $\mathcal{C}^2$);  now, we have, for $t<t_0$, $\frac{\theta(t,x_0)-\theta(t_0,x_0)}{t-t_0}\geq \frac{M(t)-M(t_0)}{t-t_0}$ so that $\displaystyle\limsup_{t\rightarrow t_0^-} \frac{M(t)-M(t_0)}{t-t_0}\leq 0$. We see that this is enough to get that $M$ is non-inecreasing on the set $\{t\ /\ M(t)>0\}$, and thus to get $M(t)\leq M(0)$; a similar argument gives $m(t)\geq m(0)$. This gives us that $\|\theta\|_\infty\leq \|\theta_0*\varphi_\eta\|_\infty\leq \|\theta_0\|_\infty$ and, if $\theta_0\geq 0$, then $\theta(x,t)\geq 0$ for all $t>0$.

$\ $ 

\noindent{\bf Third step :  global existence for the regularized problem}\\

In order to show that the $H^3$ norm of a solution   $\theta$ to equation $ (\mathcal{T}^{\epsilon,\eta}) $ does not blow up, we now compute $\partial_t ( \|  \theta\|_2^2+\|\partial_x^3\theta\|_2^2)$. As $\partial_x^3 \theta$ belongs (locally in time on $[0,T^*_{\epsilon,\eta})$) to $L^2([0,T^*_{\epsilon,\eta}) ,H^1)$ and $\partial_t \partial_x^3\theta$ to $L^2 H^{-1}$, therefore we may write
\begin{equation*}\begin{split}
\partial_t ( \|  \theta\|_2^2+\|\partial_x^3\theta\|_2^2)=& 2\int \partial_t\theta (\theta-\partial_x^6 \theta)\, dx\\ =& -2 \|\Lambda^{1/2}\theta\|_2^2-2\|\Lambda^{7/2}\theta\|_2^2-2\epsilon \|\partial_x \theta\|_2^2-2\epsilon \|\partial_x^4\theta\|_2^2
\\ & -2 \int \theta   \mathcal{H}\theta \partial_x\theta\, dx + 2 \int \partial_x^3\theta \partial_x^3(\mathcal{H}\theta \partial_x\theta)\, dx\\ =& -2 \|\Lambda^{1/2}\theta\|_2^2-2\|\Lambda^{7/2}\theta\|_2^2-2\epsilon \|\partial_x \theta\|_2^2-2\epsilon \|\partial_x^4\theta\|_2^2
\\ & -2 \int \theta   \mathcal{H}\theta \partial_x\theta\, dx + 2 \int \partial_x^3\theta \partial_x^3(\mathcal{H}\theta)\  \partial_x\theta\, dx\\ &+ 6  \int \partial_x^3\theta \partial_x^2(\mathcal{H}\theta)\  \partial_x^2\theta\, dx
+ 5  \int \partial_x^3\theta \partial_x(\mathcal{H}\theta)\  \partial_x^3\theta\, dx
\\ \leq & -2 \|\Lambda^{1/2}\theta\|_2^2-2\|\Lambda^{7/2}\theta\|_2^2-2\epsilon \|\partial_x \theta\|_2^2-2\epsilon \|\partial_x^4\theta\|_2^2
\\ & +2 \|\theta\|_\infty \|\theta\|_2\|\partial_x\theta\|_2  +(2 \|\partial_x\theta\|_7 +5\|\mathcal{H}\partial_x\theta\|_7)  \|\partial_x^3\theta\|_{7/3}^2\\ &+ 6 \|\partial_x^2\theta\|_3^2  \|\mathcal{H}\partial_x^2\theta\|_3  \end{split}\end{equation*}
We then use the boundedness of the Hilbert transform on $L^3$ and $L^7$ and the Gagliardo--Nirenberg inequalities  
$$ \|\partial_x^2 \theta\|_3\leq \|\theta\|_\infty^{1/3} \| \partial_x^3 \theta\|_2^{2/3}$$
$$ \|\partial_x\theta\|_7 \leq \|\theta\|_\infty^{5/7} \|\Lambda^{7/2}\theta\|_2^{2/7}$$
$$ \|\partial_x^3\theta\|_{7/3} \leq \|\theta\|_\infty^{1/7} \|\Lambda^{7/2}\theta\|_2^{6/7}$$
and we find, for a constant $C_0$ (that does not depend on $\theta_0$ nor on $\epsilon$),
\begin{equation}\label{h3}
\partial_t ( \|  \theta\|_2^2+\|\partial_x^3\theta\|_2^2)  \leq    C_0\|\theta_0 \|_\infty ( \|  \theta\|_2^2+\|\partial_x^3\theta\|_2^2) +2 (C_0\|\theta_0\|_\infty-1)\|\Lambda^{7/2}\theta\|_2^2 -2\epsilon \|\partial_x^4\theta\|_2^2
 \end{equation}
Thus, if $C_0 \|\theta_0\|_\infty<1$, we find that, on $[0,T^*_{\epsilon,\eta})$, we have $$ \|  \theta\|_2^2+\|\partial_x^3\theta\|_2^2 \leq e^{C_0 \|\theta_0\|_\infty t} ( \|  \theta_0*\varphi_\eta\|_2^2+\|\theta_0*\partial_x^3\varphi_\eta\|_2^2)$$
and thus $T^*_{\epsilon,\eta}=+\infty$. \\

\noindent{\bf Fourth step :  relaxing $\epsilon$}\\

From inequality \ref{h3}, we get that $\theta_{\epsilon,\eta}$ is controlled, on each bounded interval of time $[0,T]$, uniformly with respect to $\epsilon$, in the following ways :
\begin{itemize}
\item $\displaystyle\sup_{\epsilon>0} \sup_{0<t<T} \|\theta_{\epsilon,\eta}(t,.)\|_{H^3}<+\infty$
\item $\displaystyle\sup_{\epsilon>0} \int_0^T \|\theta_{\epsilon,\eta}\|_{H^{7/2}}^2\, dt <+\infty$\end{itemize} and we get from   equation $ (\mathcal{T}_{1}^{\epsilon,\eta}) $, that
\begin{itemize}
\item $\displaystyle\sup_{0<\epsilon<1} \int_0^T \|\partial_t\theta_{\epsilon,\eta}\|_{H^{1/2}}^2\, dt <+\infty$
\end{itemize}

We then use the Rellich theorem \cite{PGLR} to get that there exists a sequence $\epsilon_k\rightarrow 0$ so that $\theta_{\epsilon_k,\eta}$ converges strongly in $L^2_{\rm loc}((0,+\infty)\times\mathbb{R})$ to a limit $\theta_\eta$. As $\theta_{\epsilon,\eta} $ is (locally) bounded in $L^2 H^{7/2}$, the strong convergence
holds as well in $(L^2 H^1)_{\rm loc}$, so that $\theta_\eta$ is a solution of $ (\mathcal{T}_1) $, with initial value $\theta_0*\varphi_\eta$. \\

Moreover, we know that $\|\theta_\eta\|_\infty\leq \|\theta_0\|_\infty$ and  that, for every finite $T>0$,  
$$ \sup_{0<t<T}  \|\theta_\eta(t,.)\|_{H^3}<+\infty \text{ and } \int_0^T \|\theta_\eta\|_{H^{7/2}}^2\, dt<+\infty.$$

$\ $ 

\noindent{\bf Fifth step :  uniform estimates in  $H^{1/2}$ and $H^1$}\\

\begin{itemize}

\item control of the $L^2$ norm  : 
\begin{equation}\begin{split} \label{l2}
 \frac{1}{2}\frac{d}{dt} \left(\int \theta_\eta^2\, dx\right)= & \int \theta_ \eta\, \partial_t\theta_\eta\, dx =-\int \theta_ \eta   \Lambda \theta_ \eta  \, dx -\int \theta_ \eta  (\mathcal{H}\theta)\partial_x\theta_ \eta\, dx 
 \\ \leq &  -\int \vert \Lambda^{1/2} \theta_\eta \vert^2\, dx + \|\theta_0\|_\infty \|\theta_\eta\|_2 \|\Lambda\theta_\eta\|_2
\end{split}\end{equation}

\item  control of the $\dot H^{1/2}$ norm  :      
 \begin{equation*}\begin{split}  \frac{1}{2}\frac{d}{dt} \left(\int \vert\Lambda^{1/2} \theta_\eta\vert^2\, dx\right) =&  \int \Lambda  \theta_\eta \, \partial_t \theta_\eta\, dx\\
  =&-\int \vert  \Lambda  \theta_\eta \vert^2\, dx-\int   (\mathcal{H}\theta_\eta) \, (\Lambda\theta_\eta\,  \partial_x\theta_\eta)\, dx
  \\
  =&-\int \vert \Lambda  \theta_\eta\vert^2\, dx + \int    \theta_\eta \,  \mathcal{H}(\Lambda\theta_\eta\,  \partial_x\theta_\eta)\, dx
\end{split}
\end{equation*} 

We now use the identity, valid  for every  $f \in L^{2}$,  
\begin{equation}  
2 \mathcal{H} (f \mathcal{H}f)(x)=(\mathcal{H}f(x))^2  - f(x)^2,
\end{equation} along with,  
$$   \partial_x\theta_\eta=  \mathcal{H}\Lambda\theta_\eta, $$
to get,
$$  \| \mathcal{H}(\Lambda\theta_\eta\,  \partial_x\theta_\eta)\|_1\leq  \|\Lambda\theta_\eta\|_2^2,$$
  and finally obtain
  \begin{equation}\label{h1/2}  \frac{d}{dt} \int \vert \Lambda^{1/2}\theta_\eta \vert ^2\, dx+2(1- \|\theta_0\|_\infty) \int   \vert \Lambda  \theta_\eta \vert^2\, dx  \leq 0. \end{equation}    

\item control of the $\dot H^{1}$ norm  :  we  write    
 \begin{equation*}\begin{split} 
  \frac{1}{2}\frac{d}{dt} \int \vert \Lambda  \theta_\eta \vert^2\, dx =&  \int \Lambda^2  \theta_\eta \, \partial_t \theta_\eta\, dx\\
  =&-\int  \vert  \Lambda^{3/2}  \theta_\eta \vert^2\, dx- \frac{1}{2}\int  \partial_x (\mathcal{H}\theta_\eta) \, (\partial_x\theta_\eta)^2\, dx.
\end{split}
\end{equation*}
Using a Gagliardo--Nirenberg inequality, we get
$$ \frac{1}{2} \left\vert \int  \partial_x (\mathcal{H}\theta_\eta) \, (\partial_x\theta_\eta)^2\, dx\right\vert\leq C \|\partial_x\theta\|_3^3\leq C_1 \|\theta\|_\infty \|\Lambda^{3/2}\theta_\eta\|_2^2,$$ and finally obtain,
  \begin{equation} \label{sob}  \frac{d}{dt} \left(\int \vert\Lambda \theta_\eta \vert^2\, dx\right)+2(1-C_1\|\theta_0\|_\infty) \int  \vert \Lambda ^{3/2}  \theta_\eta\vert^2\, dx  \leq 0. \end{equation} 
\end{itemize}

$\ $ 

\noindent{\bf Sixth step :  relaxing $\eta$}\\

From inequalities (\ref{l2}) and  (\ref{h1/2}), we get that, for $\theta_0\in H^{1/2}$, (when $\|\theta_0\|_\infty$ is small enough)  $\theta_{\eta}$ is controlled, on each bounded interval of time $[0,T]$, uniformly with respect to $\eta$, in the following ways :
\begin{itemize}
\item $\displaystyle\sup_{\eta>0} \sup_{0<t<T} \|\theta_{\eta}(t,.)\|_{H^{1/2}}<+\infty$,
\item $\displaystyle\sup_{\eta>0} \int_0^T \|\theta_{\eta}\|_{H^{1}}^2\, dt <+\infty$\end{itemize} and we get from   equation $ (\mathcal{T}_1)$, that
\begin{itemize}
\item $\displaystyle\sup_{\eta>0} \int_0^T \|\partial_t\theta_{\epsilon,\eta}\|_{H^{-1/4}}^2\, dt <+\infty$.
\end{itemize}

We may then use the Rellich theorem  \cite{PGLR} and  get that there exists a sequence $\eta_k\rightarrow 0$ so that $\theta_{\eta_k}$ converges strongly in $L^2_{\rm loc}((0,+\infty)\times\mathbb{R})$ to a limit $\theta$. As $\theta_{\eta} $ is (locally) bounded in $L^2 H^{1}$,we have weak convergence in $L^2 H^1$; we then write $\mathcal{H}\theta_\eta \partial_x\theta_\eta=\partial_x(\theta_\eta \mathcal{H}\theta_\eta)-\theta_\eta \mathcal{H}\partial_x\theta_\eta$ and find that $\theta$ is a solution of $ (\mathcal{T}_1) $, with initial value $\theta_0$. \\

Moreover, we find  that we have 
\begin{itemize}
\item  $\|\theta \|_\infty\leq \|\theta_0\|_\infty$
\item $\displaystyle\sup_{t>0} \| \Lambda^{1/2} \theta(t,.)\|_2 \leq \|  \Lambda^{1/2} \theta_0\|_2$
\item $\int_0^{+\infty} \|\Lambda \theta\|_2^2\, ds \leq \frac{1}{2(1- \|\theta_0\|_\infty)} \|\Lambda^{1/2}\theta_0\|_2^2 $
\item $\|\theta(t,.)\|_2\leq \| \theta_0\|_2 + \|\theta_0\|_\infty \int_0^t \|\Lambda \theta(s,.)\|_2\, ds.$
\end{itemize}

Similarly, if  $\theta_0\in H^{1}$  (with  $\|\theta_0\|_\infty$   small enough) , then inequality (\ref{sob})  will give a control of the $H^1$ norm of $\theta_\eta$ uniformly with respect to $\eta$, and thus, we find for the limit $\theta$  that,
\begin{itemize} 
\item $\displaystyle\sup_{t>0} \| \Lambda \theta(t,.)\|_2 \leq \|  \Lambda  \theta_0\|_2$,
\item $\int_0^{+\infty} \|\Lambda^{3/2} \theta\|_2^2\, ds \leq \frac{1}{2(1-C_1 \|\theta_0\|_\infty)} \|\Lambda \theta_0\|_2^2$. \end{itemize}

 \vskip0.2cm\noindent{\bf Acknowledgment}:  The first author thanks Diego C\'ordoba for useful discussions regarding this model. He was partially supported by the ERC grant Stg-203138-CDSIF  and the National Grant MTM2014-59488-P from the Spanish government.  The authors which to thank the anonymous referees for their careful reading of the article.

\bibliographystyle{amsplain}

\begin{thebibliography}{10}

\bibitem[1]{Bak} G.R. Baker, X. Li,  A.C. Morlet. \emph{Analytic structure of
two 1D-transport equations with nonlocal fluxes}.
Physica D: Nonlinear Phenomena, {\bf{91}}(4):349-375, 1996.

\bibitem[2]{BG}	H. Bae, R. Granero-Belinch\'on. \emph{Global existence for some transport equations with nonlocal velocity.}	Advances in Mathematics, vol. {\bf{269}}, pp 197-219, 2015.


\bibitem[3]{CV} L.A. Caffarelli, A. Vasseur. \emph{Drift diffusion equations with fractional diffusion and the quasi-geostrophic equation}. Annals of Math. (2), {\bf{171}}(3):1903-1930, 2010.

\bibitem[4]{CaC} A. Castro, D. C\'ordoba. \emph{Global existence, singularities and Ill-posedness for a non-local flux.}   Advances in Mathematics. {\bf{219}} (2008), 6, 1916-1936.

\bibitem[5]{CaC2} A. Castro, D. C\'ordoba. \emph{Infinite energy solutions of the surface quasi-geostrophic equation.} Advances in Mathematics. {\bf{225}} (2010) 1820-1829.


\bibitem[6]{ChCCF} D. Chae, A. C\'ordoba, D. C\'ordoba, M. A. Fontelos. \emph{Finite time singularities in a 
1D model of the quasi-geostrophic equation.} Advances in Mathematics. {\bf{194}} (2005), 203-223.

\bibitem[7]{CCS} C. H. Chan, M. Czubac, L. Silvestre. \emph{Eventual regularization of the slightly supercritical fractional Burgers equation. Discrete and Continuous Dynamical Systems}, Volume:  {\bf{27}}, Number: {{2}}, June 2010, Pages 847-861.



\bibitem[8] {CM2}   R. Coifman, Y. Meyer.  \emph {Wavelets: Calder\'on-Zygmund and Multilinear Operators}, Cambridge University Press, 336 pages

\bibitem[9]{CLM}  P. Constantin, P. Lax, A. Majda. \textit{A simple one-dimensional
model for the three dimensional vorticity}, Comm.\ Pure
Appl.\ Math.\ {\bf{38}} (1985), 715-724.


\bibitem[10]{CMT} P. Constantin, A.J.  Majda, E. Tabak. \textit{Formation of strong fronts in
  the $2$-D quasigeostrophic thermal active scalar}.
  {\newblock}Nonlinearity, {\bf{7}} (1994), pp. 1495-1533.


\bibitem[11]{CVi} P. Constantin, V. Vicol. \emph{Nonlinear maximum principles for dissipative linear nonlocal operators and applications}. Geometric And Functional Analysis, {\bf{22}}(5):1289-1321, 2012.

 \bibitem[12]{CC}A. C\'ordoba, D. C\'ordoba.
  \emph{A maximum principle applied to quasi-geostrophic equations,}
  {\newblock}Comm. Math. Phys. {\bf{249}} (2004), pp. 511-528.

\bibitem[13]{CCF2} A. C\'ordoba, D. C\'ordoba, M.A. Fontelos. \emph{Integral inequalities for the Hilbert transform applied to a nonlocal transport equation}. J. Math. Pures Appl. (9), {\bf{86}}, 6:529-540, 2006



\bibitem[14]{CCF} A. C\'ordoba, D. C\'ordoba,  M.A. Fontelos. \textit{Formation of singularities 
for a transport equation with nonlocal velocity}, Annals of Math. 162 (2005) 
{\bf{3}}, 1375-1387.

\bibitem[15]{De1} S. De Gregorio. \textit{On a one-dimensional model for the
three-dimensional vorticity equation}, J. Statist.\ Phys.
{\bf{59}} (1990), 1251-1263.

\bibitem[16]{Do} T. Do. \emph{On a 1d transport equation with nonlocal velocity and supercritical dissipation}. Journal of Differential Equations, {\bf{256}}-9:3166-3178, 2014.


 \bibitem[17]{Dong}    H. Dong. \emph{Well-posedness for a transport equation with 
nonlocal velocity}. J. Funct. Anal, {\bf{255}}:3070-3097, (2008).


 \bibitem[18]{Hed} L. Hedberg. \emph{On certain convolution inequalities}. Proc. Amer. Math. Soc. {\bf{10}} (1972), 505-510.


\bibitem[19]{HMW}   R. Hunt, B. Muckenhoupt,  R. Wheeden. \emph{Weighted norm inequalities for 
the conjugate function and Hilbert transform}, Trans. Amer. Math. Soc. {\bf{176}} (1973), 227-251.


\bibitem[20]{Ki}  A. Kiselev. \emph{Regularity and blow up for active scalars}.
Math. Model. Nat. Phenom, 5(4):225-255, 2010

\bibitem[21]{K} A. Kiselev. \emph{ Nonlocal maximum principles for active scalars}. Advances in Mathematics, 227,{\bf{5}}:1806-1826, 2011.

\bibitem[22]{KNS}  A. Kiselev, F. Nazarov and R. Shterenberg. \emph{On blow up and regularity in dissipative Burgers equation}, Dynamics of PDEs, {\bf{5}} (2008), 211-240

\bibitem[23]{KNV}  A. Kiselev, F. Nazarov, and A. Volberg. \emph{Global well-posedness for the critical 2D dissipative quasi-geostrophic equation}.
Invent. Math., {\bf{167}}(3):445-453, 2007

\bibitem[24] {Laz} O. Lazar. \emph{On a 1D nonlocal transport equation with nonlocal velocity and subcritical or supercritical diffusion}, preprint.


\bibitem[25] {PGLR} P. G. Lemari\'e-Rieusset. \emph{Recent developments in the Navier-Stokes problem.} Chapman \& Hall/CRC (2002).




\bibitem[26]{RL}  D. Li, J. L. Rodrigo. {Blow-up of solutions for a 1D transport equation with nonlocal velocity and supercritical dissipation.} {\it Advances in Mathematics}, {\bf{217}}, no. 6, 2563-2568 (2008).

 
\bibitem[27]{RL2} D. Li, J. L. Rodrigo. \emph{On a One-Dimensional Nonlocal Flux with Fractional Dissipation}, SIAM J. Math. Anal. {\bf{43}} (2011), 507-526.  

\bibitem[28]{Muc} B. Muckenhoupt.  \emph{Weighted norm inequalities for the Hardy maximal function}. Transactions of the American Mathematical Society, vol. {\bf{165}}: 207-226. (1972).

\bibitem[29]{OSW} H. Okamoto, T. Sakajo, M. Wunsch.  {On a generalization of the Constantin-Lax-Majda equation.} Nonlinearity, {\bf{21}}(10): 2447-2461 (2008).

\bibitem[30]{Stein}  E. Stein. {Harmonic Analysis : Real Variable Methods, Orthogonality and Oscillatory Integrals}, Princeton Math. Series {\bf{43}}, Princeton Univ. Press, Princeton, NJ, 1993.


\bibitem[31]{S} L. Silvestre. Eventual regularization for the slightly supercritical quasi-geostrophic equation. Ann. Inst. H. Poincar\'e Anal. Non
Lin\'eaire, {\bf{27}}(2):693-704, 2010

\bibitem[32]{ViS} L. Silvestre, V. Vicol.  Transactions of the American Mathematical Society {\bf{368}} (2016), no. 9, 6159-6188..


 

  
  
  
	


















\end{thebibliography}

\vspace{1cm}

\begin{quote}
\begin{tabular}{ll}
\vspace{0,2cm}
Omar Lazar & Pierre-Gilles Lemari\'e-Rieusset\\
{\small Instituto de Ciencias Matem\'aticas (ICMAT)} & {\small Universit\'e d'Evry Val d'Essonne}\\
{\small Consejo Superior de Investigaciones Cient\'ificas} & {\small LaMME  (UMR  CNRS 8071)}\\
{\small C/ Nicolas Cabrera 13-15, 28049 Madrid, Spain} & {\small 23 Boulevard de France, 91037 \'Evry Cedex, France }\\
{\small Email: omar.lazar@icmat.es} & {\small Email: plemarie@univ-evry.fr}
\end{tabular}
\end{quote}
\end{document}